
\ifx\shlhetal\undefinedcontrolsequence\let\shlhetal\relax\fi
\def\fmtname{AmS-TeX}

\def\fmtversion{2.1}
\catcode`\@=11
\ifx\amstexloaded@\relax\catcode`\@=\active
  \endinput\else\let\amstexloaded@\relax\fi
\newlinechar=`\^^J
\def\W@{\immediate\write\sixt@@n}
\def\CR@{\W@{^^J\fmtname - Version \fmtversion^^J%
COPYRIGHT 1985, 1990, 1991 - AMERICAN MATHEMATICAL SOCIETY^^J%
Use of this macro package is not restricted provided^^J%
each use is acknowledged upon publication.^^J}}
\CR@ \everyjob{\CR@}
\message{Loading definitions for}
\message{misc utility macros,}
\toksdef\toks@@=2
\long\def\rightappend@#1\to#2{\toks@{\\{#1}}\toks@@
 =\expandafter{#2}\xdef#2{\the\toks@@\the\toks@}\toks@{}\toks@@{}}
\def\alloclist@{}
\newif\ifalloc@
\def\showallocations{{\def\\{\immediate\write\m@ne}\alloclist@}\alloc@true}
\def\alloc@#1#2#3#4#5{\global\advance\count1#1by\@ne
 \ch@ck#1#4#2\allocationnumber=\count1#1
 \global#3#5=\allocationnumber
 \edef\next@{\string#5=\string#2\the\allocationnumber}%
 \expandafter\rightappend@\next@\to\alloclist@}
\newcount\count@@
\newcount\count@@@
\def\FN@{\futurelet\next}
\def\DN@{\def\next@}
\def\DNii@{\def\nextii@}
\def\RIfM@{\relax\ifmmode}
\def\RIfMIfI@{\relax\ifmmode\ifinner}
\def\setboxz@h{\setbox\z@\hbox}
\def\wdz@{\wd\z@}
\def\boxz@{\box\z@}
\def\setbox@ne{\setbox\@ne}
\def\wd@ne{\wd\@ne}
\def\iterate{\body\expandafter\iterate\else\fi}
\def\err@#1{\errmessage{AmS-TeX error: #1}}
\newhelp\defaulthelp@{Sorry, I already gave what help I could...^^J
Maybe you should try asking a human?^^J
An error might have occurred before I noticed any problems.^^J
``If all else fails, read the instructions.''}
\def\Err@{\errhelp\defaulthelp@\err@}
\def\eat@#1{}
\def\in@#1#2{\def\in@@##1#1##2##3\in@@{\ifx\in@##2\in@false\else\in@true\fi}%
 \in@@#2#1\in@\in@@}
\newif\ifin@
\def\space@.{\futurelet\space@\relax}
\space@. %
\newhelp\athelp@
{Only certain combinations beginning with @ make sense to me.^^J
Perhaps you wanted \string\@\space for a printed @?^^J
I've ignored the character or group after @.}
{\catcode`\~=\active 
 \lccode`\~=`\@ \lowercase{\gdef~{\FN@\at@}}}
\def\at@{\let\next@\at@@
 \ifcat\noexpand\next a\else\ifcat\noexpand\next0\else
 \ifcat\noexpand\next\relax\else
   \let\next\at@@@\fi\fi\fi
 \next@}
\def\at@@#1{\expandafter
 \ifx\csname\space @\string#1\endcsname\relax
  \expandafter\at@@@ \else
  \csname\space @\string#1\expandafter\endcsname\fi}
\def\at@@@#1{\errhelp\athelp@ \err@{\Invalid@@ @}}
\def\atdef@#1{\expandafter\def\csname\space @\string#1\endcsname}
\newhelp\defahelp@{If you typed \string\define\space cs instead of
\string\define\string\cs\space^^J
I've substituted an inaccessible control sequence so that your^^J
definition will be completed without mixing me up too badly.^^J
If you typed \string\define{\string\cs} the inaccessible control sequence^^J
was defined to be \string\cs, and the rest of your^^J
definition appears as input.}
\newhelp\defbhelp@{I've ignored your definition, because it might^^J
conflict with other uses that are important to me.}
\def\define{\FN@\define@}
\def\define@{\ifcat\noexpand\next\relax
 \expandafter\define@@\else\errhelp\defahelp@                               
 \err@{\string\define\space must be followed by a control
 sequence}\expandafter\def\expandafter\nextii@\fi}                          
\def\undefined@@@@@@@@@@{}
\def\preloaded@@@@@@@@@@{}
\def\next@@@@@@@@@@{}
\def\define@@#1{\ifx#1\relax\errhelp\defbhelp@                              
 \err@{\string#1\space is already defined}\DN@{\DNii@}\else
 \expandafter\ifx\csname\expandafter\eat@\string                            
 #1@@@@@@@@@@\endcsname\undefined@@@@@@@@@@\errhelp\defbhelp@
 \err@{\string#1\space can't be defined}\DN@{\DNii@}\else
 \expandafter\ifx\csname\expandafter\eat@\string#1\endcsname\relax          
 \global\let#1\undefined\DN@{\def#1}\else\errhelp\defbhelp@
 \err@{\string#1\space is already defined}\DN@{\DNii@}\fi
 \fi\fi\next@}

\def\predefine#1#2{\let#1#2}
\def\undefine#1{\let#1\undefined}
\message{page layout,}
\newdimen\captionwidth@
\captionwidth@\hsize
\advance\captionwidth@-1.5in
\def\pagewidth#1{\hsize#1\relax
 \captionwidth@\hsize\advance\captionwidth@-1.5in}
\def\pageheight#1{\vsize#1\relax}
\def\hcorrection#1{\advance\hoffset#1\relax}
\def\vcorrection#1{\advance\voffset#1\relax}
\message{accents/punctuation,}

\let\graveaccent\`
\let\acuteaccent\'
\let\tildeaccent\~
\let\hataccent\^
\let\underscore\_
\let\B\=
\let\D\.
\let\ic@\/
\def\/{\unskip\ic@}
\def\textfonti{\the\textfont\@ne}
\def\t#1#2{{\edef\next@{\the\font}\textfonti\accent"7F \next@#1#2}}
\def~{\unskip\nobreak\ \ignorespaces}
\def\.{.\spacefactor\@m}
\atdef@;{\leavevmode\null;}
\atdef@:{\leavevmode\null:}
\atdef@?{\leavevmode\null?}
\edef\@{\string @}
\def\relaxnext@{\let\next\relax}
\atdef@-{\relaxnext@\leavevmode
 \DN@{\ifx\next-\DN@-{\FN@\nextii@}\else
  \DN@{\leavevmode\hbox{-}}\fi\next@}%
 \DNii@{\ifx\next-\DN@-{\leavevmode\hbox{---}}\else
  \DN@{\leavevmode\hbox{--}}\fi\next@}%
 \FN@\next@}
\def\srdr@{\kern.16667em}
\def\drsr@{\kern.02778em}
\def\sldl@{\drsr@}
\def\dlsl@{\srdr@}
\atdef@"{\unskip\relaxnext@
 \DN@{\ifx\next\space@\DN@. {\FN@\nextii@}\else
  \DN@.{\FN@\nextii@}\fi\next@.}%
 \DNii@{\ifx\next`\DN@`{\FN@\nextiii@}\else
  \ifx\next\lq\DN@\lq{\FN@\nextiii@}\else
  \DN@####1{\FN@\nextiv@}\fi\fi\next@}%
 \def\nextiii@{\ifx\next`\DN@`{\sldl@``}\else\ifx\next\lq
  \DN@\lq{\sldl@``}\else\DN@{\dlsl@`}\fi\fi\next@}%
 \def\nextiv@{\ifx\next'\DN@'{\srdr@''}\else
  \ifx\next\rq\DN@\rq{\srdr@''}\else\DN@{\drsr@'}\fi\fi\next@}%
 \FN@\next@}

\def\textfontii{\the\textfont\tw@}
\def\lbrace@{\delimiter"4266308 }
\def\rbrace@{\delimiter"5267309 }
\def\{{\RIfM@\lbrace@\else{\textfontii f}\spacefactor\@m\fi}
\def\}{\RIfM@\rbrace@\else
 \let\@sf\empty\ifhmode\edef\@sf{\spacefactor\the\spacefactor}\fi
 {\textfontii g}\@sf\relax\fi}
\let\lbrace\{
\let\rbrace\}
\def\AmSTeX{{\textfontii A\kern-.1667em%
  \lower.5ex\hbox{M}\kern-.125emS}-\TeX}
\message{line and page breaks,}
\def\vmodeerr@#1{\Err@{\string#1\space not allowed between paragraphs}}
\def\mathmodeerr@#1{\Err@{\string#1\space not allowed in math mode}}
\def\linebreak{\RIfM@\mathmodeerr@\linebreak\else
 \ifhmode\unskip\unkern\break\else\vmodeerr@\linebreak\fi\fi}

\newskip\saveskip@
\def\allowlinebreak{\RIfM@\mathmodeerr@\allowlinebreak\else
 \ifhmode\saveskip@\lastskip\unskip
 \allowbreak\ifdim\saveskip@>\z@\hskip\saveskip@\fi
 \else\vmodeerr@\allowlinebreak\fi\fi}
\def\nolinebreak{\RIfM@\mathmodeerr@\nolinebreak\else
 \ifhmode\saveskip@\lastskip\unskip
 \nobreak\ifdim\saveskip@>\z@\hskip\saveskip@\fi
 \else\vmodeerr@\nolinebreak\fi\fi}
\def\newline{\relaxnext@
 \DN@{\RIfM@\expandafter\mathmodeerr@\expandafter\newline\else
  \ifhmode\ifx\next\par\else
  \expandafter\unskip\expandafter\null\expandafter\hfill\expandafter\break\fi
  \else
  \expandafter\vmodeerr@\expandafter\newline\fi\fi}%
 \FN@\next@}
\def\dmatherr@#1{\Err@{\string#1\space not allowed in display math mode}}
\def\nondmatherr@#1{\Err@{\string#1\space not allowed in non-display math
 mode}}
\def\onlydmatherr@#1{\Err@{\string#1\space allowed only in display math mode}}
\def\nonmatherr@#1{\Err@{\string#1\space allowed only in math mode}}
\def\mathbreak{\RIfMIfI@\break\else
 \dmatherr@\mathbreak\fi\else\nonmatherr@\mathbreak\fi}
\def\nomathbreak{\RIfMIfI@\nobreak\else
 \dmatherr@\nomathbreak\fi\else\nonmatherr@\nomathbreak\fi}
\def\allowmathbreak{\RIfMIfI@\allowbreak\else
 \dmatherr@\allowmathbreak\fi\else\nonmatherr@\allowmathbreak\fi}
\def\pagebreak{\RIfM@
 \ifinner\nondmatherr@\pagebreak\else\postdisplaypenalty-\@M\fi
 \else\ifvmode\removelastskip\break\else\vadjust{\break}\fi\fi}
\def\nopagebreak{\RIfM@
 \ifinner\nondmatherr@\nopagebreak\else\postdisplaypenalty\@M\fi
 \else\ifvmode\nobreak\else\vadjust{\nobreak}\fi\fi}
\def\nonvmodeerr@#1{\Err@{\string#1\space not allowed within a paragraph
 or in math}}
\def\vnonvmode@#1#2{\relaxnext@\DNii@{\ifx\next\par\DN@{#1}\else
 \DN@{#2}\fi\next@}%
 \ifvmode\DN@{#1}\else
 \DN@{\FN@\nextii@}\fi\next@}
\def\newpage{\vnonvmode@{\vfill\break}{\nonvmodeerr@\newpage}}
\def\smallpagebreak{\vnonvmode@\smallbreak{\nonvmodeerr@\smallpagebreak}}
\def\medpagebreak{\vnonvmode@\medbreak{\nonvmodeerr@\medpagebreak}}
\def\bigpagebreak{\vnonvmode@\bigbreak{\nonvmodeerr@\bigpagebreak}}
\def\NoBlackBoxes{\global\overfullrule\z@}
\def\BlackBoxes{\global\overfullrule5\p@}
\def\Invalid@#1{\def#1{\Err@{\Invalid@@\string#1}}}
\def\Invalid@@{Invalid use of }
\message{figures,}
\Invalid@\caption
\Invalid@\captionwidth
\newdimen\smallcaptionwidth@
\def\topspace{\mid@false\ins@}
\def\midspace{\mid@true\ins@}
\newif\ifmid@
\def\captionfont@{}
\def\ins@#1{\relaxnext@\allowbreak
 \smallcaptionwidth@\captionwidth@\gdef\thespace@{#1}%
 \DN@{\ifx\next\space@\DN@. {\FN@\nextii@}\else
  \DN@.{\FN@\nextii@}\fi\next@.}%
 \DNii@{\ifx\next\caption\DN@\caption{\FN@\nextiii@}%
  \else\let\next@\nextiv@\fi\next@}%
 \def\nextiv@{\vnonvmode@
  {\ifmid@\expandafter\midinsert\else\expandafter\topinsert\fi
   \vbox to\thespace@{}\endinsert}
  {\ifmid@\nonvmodeerr@\midspace\else\nonvmodeerr@\topspace\fi}}%
 \def\nextiii@{\ifx\next\captionwidth\expandafter\nextv@
  \else\expandafter\nextvi@\fi}%
 \def\nextv@\captionwidth##1##2{\smallcaptionwidth@##1\relax\nextvi@{##2}}%
 \def\nextvi@##1{\def\thecaption@{\captionfont@##1}%
  \DN@{\ifx\next\space@\DN@. {\FN@\nextvii@}\else
   \DN@.{\FN@\nextvii@}\fi\next@.}%
  \FN@\next@}%
 \def\nextvii@{\vnonvmode@
  {\ifmid@\expandafter\midinsert\else
  \expandafter\topinsert\fi\vbox to\thespace@{}\nobreak\smallskip
  \setboxz@h{\noindent\ignorespaces\thecaption@\unskip}%
  \ifdim\wdz@>\smallcaptionwidth@\centerline{\vbox{\hsize\smallcaptionwidth@
   \noindent\ignorespaces\thecaption@\unskip}}%
  \else\centerline{\boxz@}\fi\endinsert}
  {\ifmid@\nonvmodeerr@\midspace
  \else\nonvmodeerr@\topspace\fi}}%
 \FN@\next@}
\message{comments,}
\def\newcodes@{\catcode`\\12\catcode`\{12\catcode`\}12\catcode`\#12%
 \catcode`\%12\relax}
\def\oldcodes@{\catcode`\\0\catcode`\{1\catcode`\}2\catcode`\#6%
 \catcode`\%14\relax}
\def\comment{\newcodes@\endlinechar=10 \comment@}
{\lccode`\0=`\\
\lowercase{\gdef\comment@#1^^J{\comment@@#10endcomment\comment@@@}%
\gdef\comment@@#10endcomment{\FN@\comment@@@}%
\gdef\comment@@@#1\comment@@@{\ifx\next\comment@@@\let\next\comment@
 \else\def\next{\oldcodes@\endlinechar=`\^^M\relax}%
 \fi\next}}}
\def\pr@m@s{\ifx'\next\DN@##1{\prim@s}\else\let\next@\egroup\fi\next@}
\def\prime{{\null\prime@\null}}
\mathchardef\prime@="0230
\let\dsize\displaystyle

\let\ssize\scriptstyle

\message{math spacing,}
\def\,{\RIfM@\mskip\thinmuskip\relax\else\kern.16667em\fi}
\def\!{\RIfM@\mskip-\thinmuskip\relax\else\kern-.16667em\fi}
\let\thinspace\,
\let\negthinspace\!
\def\medspace{\RIfM@\mskip\medmuskip\relax\else\kern.222222em\fi}
\def\negmedspace{\RIfM@\mskip-\medmuskip\relax\else\kern-.222222em\fi}
\def\thickspace{\RIfM@\mskip\thickmuskip\relax\else\kern.27777em\fi}
\let\;\thickspace
\def\negthickspace{\RIfM@\mskip-\thickmuskip\relax\else
 \kern-.27777em\fi}
\atdef@,{\RIfM@\mskip.1\thinmuskip\else\leavevmode\null,\fi}
\atdef@!{\RIfM@\mskip-.1\thinmuskip\else\leavevmode\null!\fi}
\atdef@.{\RIfM@&&\else\leavevmode.\spacefactor3000 \fi}
\def\and{\DOTSB\;\mathchar"3026 \;}

\message{fractions,}
\def\frac#1#2{{#1\over#2}}

\newdimen\ex@
\ex@.2326ex
\Invalid@\thickness
\def\thickfrac{\relaxnext@
 \DN@{\ifx\next\thickness\let\next@\nextii@\else
 \DN@{\nextii@\thickness1}\fi\next@}%
 \DNii@\thickness##1##2##3{{##2\above##1\ex@##3}}%
 \FN@\next@}

\def\thickfracwithdelims#1#2{\relaxnext@\def\ldelim@{#1}\def\rdelim@{#2}%
 \DN@{\ifx\next\thickness\let\next@\nextii@\else
 \DN@{\nextii@\thickness1}\fi\next@}%
 \DNii@\thickness##1##2##3{{##2\abovewithdelims
 \ldelim@\rdelim@##1\ex@##3}}%
 \FN@\next@}
\def\binom#1#2{{#1\choose#2}}

\def\:{\nobreak\hskip.1111em\mathpunct{}\nonscript\mkern-\thinmuskip{:}\hskip
 .3333emplus.0555em\relax}
\def\snug{\unskip\kern-\mathsurround}
\message{smash commands,}
\def\topsmash{\top@true\bot@false\smash@}
\def\botsmash{\top@false\bot@true\smash@}
\newif\iftop@
\newif\ifbot@
\def\smash{\top@true\bot@true\smash@}
\def\smash@{\RIfM@\expandafter\mathpalette\expandafter\mathsm@sh\else
 \expandafter\makesm@sh\fi}
\def\finsm@sh{\iftop@\ht\z@\z@\fi\ifbot@\dp\z@\z@\fi\leavevmode\boxz@}
\message{large operator symbols,}
\def\LimitsOnSums{\global\let\slimits@\displaylimits}
\def\NoLimitsOnSums{\global\let\slimits@\nolimits}
\LimitsOnSums
\mathchardef\coprod@="1360       \def\coprod{\DOTSB\coprod@\slimits@}
\mathchardef\bigvee@="1357       \def\bigvee{\DOTSB\bigvee@\slimits@}
\mathchardef\bigwedge@="1356     \def\bigwedge{\DOTSB\bigwedge@\slimits@}
\mathchardef\biguplus@="1355     \def\biguplus{\DOTSB\biguplus@\slimits@}
\mathchardef\bigcap@="1354       \def\bigcap{\DOTSB\bigcap@\slimits@}
\mathchardef\bigcup@="1353       \def\bigcup{\DOTSB\bigcup@\slimits@}
\mathchardef\prod@="1351         \def\prod{\DOTSB\prod@\slimits@}
\mathchardef\sum@="1350          \def\sum{\DOTSB\sum@\slimits@}
\mathchardef\bigotimes@="134E    \def\bigotimes{\DOTSB\bigotimes@\slimits@}
\mathchardef\bigoplus@="134C     \def\bigoplus{\DOTSB\bigoplus@\slimits@}
\mathchardef\bigodot@="134A      \def\bigodot{\DOTSB\bigodot@\slimits@}
\mathchardef\bigsqcup@="1346     \def\bigsqcup{\DOTSB\bigsqcup@\slimits@}
\message{integrals,}
\def\LimitsOnInts{\global\let\ilimits@\displaylimits}
\def\NoLimitsOnInts{\global\let\ilimits@\nolimits}
\NoLimitsOnInts
\def\int{\DOTSI\intop\ilimits@}
\def\oint{\DOTSI\ointop\ilimits@}
\def\intic@{\mathchoice{\hskip.5em}{\hskip.4em}{\hskip.4em}{\hskip.4em}}
\def\negintic@{\mathchoice
 {\hskip-.5em}{\hskip-.4em}{\hskip-.4em}{\hskip-.4em}}
\def\intkern@{\mathchoice{\!\!\!}{\!\!}{\!\!}{\!\!}}
\def\intdots@{\mathchoice{\plaincdots@}
 {{\cdotp}\mkern1.5mu{\cdotp}\mkern1.5mu{\cdotp}}
 {{\cdotp}\mkern1mu{\cdotp}\mkern1mu{\cdotp}}
 {{\cdotp}\mkern1mu{\cdotp}\mkern1mu{\cdotp}}}
\newcount\intno@
\def\iint{\DOTSI\intno@\tw@\FN@\ints@}
\def\iiint{\DOTSI\intno@\thr@@\FN@\ints@}
\def\iiiint{\DOTSI\intno@4 \FN@\ints@}
\def\idotsint{\DOTSI\intno@\z@\FN@\ints@}
\def\ints@{\findlimits@\ints@@}
\newif\iflimtoken@
\newif\iflimits@
\def\findlimits@{\limtoken@true\ifx\next\limits\limits@true
 \else\ifx\next\nolimits\limits@false\else
 \limtoken@false\ifx\ilimits@\nolimits\limits@false\else
 \ifinner\limits@false\else\limits@true\fi\fi\fi\fi}
\def\multint@{\int\ifnum\intno@=\z@\intdots@                                
 \else\intkern@\fi                                                          
 \ifnum\intno@>\tw@\int\intkern@\fi                                         
 \ifnum\intno@>\thr@@\int\intkern@\fi                                       
 \int}                                                                      
\def\multintlimits@{\intop\ifnum\intno@=\z@\intdots@\else\intkern@\fi
 \ifnum\intno@>\tw@\intop\intkern@\fi
 \ifnum\intno@>\thr@@\intop\intkern@\fi\intop}
\def\ints@@{\iflimtoken@                                                    
 \def\ints@@@{\iflimits@\negintic@\mathop{\intic@\multintlimits@}\limits    
  \else\multint@\nolimits\fi                                                
  \eat@}                                                                    
 \else                                                                      
 \def\ints@@@{\iflimits@\negintic@
  \mathop{\intic@\multintlimits@}\limits\else
  \multint@\nolimits\fi}\fi\ints@@@}
\def\LimitsOnNames{\global\let\nlimits@\displaylimits}
\def\NoLimitsOnNames{\global\let\nlimits@\nolimits@}
\LimitsOnNames
\def\nolimits@{\relaxnext@
 \DN@{\ifx\next\limits\DN@\limits{\nolimits}\else
  \let\next@\nolimits\fi\next@}%
 \FN@\next@}
\message{operator names,}
\def\newmcodes@{\mathcode`\'"27\mathcode`\*"2A\mathcode`\."613A%
 \mathcode`\-"2D\mathcode`\/"2F\mathcode`\:"603A }
\def\operatorname#1{\mathop{\newmcodes@\kern\z@\fam\z@#1}\nolimits@}
\def\operatornamewithlimits#1{\mathop{\newmcodes@\kern\z@\fam\z@#1}\nlimits@}
\def\qopname@#1{\mathop{\fam\z@#1}\nolimits@}
\def\qopnamewl@#1{\mathop{\fam\z@#1}\nlimits@}
\def\arccos{\qopname@{arccos}}
\def\arcsin{\qopname@{arcsin}}
\def\arctan{\qopname@{arctan}}
\def\arg{\qopname@{arg}}
\def\cos{\qopname@{cos}}
\def\cosh{\qopname@{cosh}}
\def\cot{\qopname@{cot}}
\def\coth{\qopname@{coth}}
\def\csc{\qopname@{csc}}
\def\deg{\qopname@{deg}}
\def\det{\qopnamewl@{det}}
\def\dim{\qopname@{dim}}
\def\exp{\qopname@{exp}}
\def\gcd{\qopnamewl@{gcd}}
\def\hom{\qopname@{hom}}
\def\inf{\qopnamewl@{inf}}
\def\injlim{\qopnamewl@{inj\,lim}}
\def\ker{\qopname@{ker}}
\def\lg{\qopname@{lg}}
\def\lim{\qopnamewl@{lim}}
\def\liminf{\qopnamewl@{lim\,inf}}
\def\limsup{\qopnamewl@{lim\,sup}}
\def\ln{\qopname@{ln}}
\def\log{\qopname@{log}}
\def\max{\qopnamewl@{max}}
\def\min{\qopnamewl@{min}}
\def\Pr{\qopnamewl@{Pr}}
\def\projlim{\qopnamewl@{proj\,lim}}
\def\sec{\qopname@{sec}}
\def\sin{\qopname@{sin}}
\def\sinh{\qopname@{sinh}}
\def\sup{\qopnamewl@{sup}}
\def\tan{\qopname@{tan}}
\def\tanh{\qopname@{tanh}}
\def\varinjlim{\mathop{\vtop{\ialign{##\crcr
 \hfil\rm lim\hfil\crcr\noalign{\nointerlineskip}\rightarrowfill\crcr
 \noalign{\nointerlineskip\kern-\ex@}\crcr}}}}
\def\varprojlim{\mathop{\vtop{\ialign{##\crcr
 \hfil\rm lim\hfil\crcr\noalign{\nointerlineskip}\leftarrowfill\crcr
 \noalign{\nointerlineskip\kern-\ex@}\crcr}}}}
\def\varliminf{\mathop{\underline{\vrule height\z@ depth.2exwidth\z@
 \hbox{\rm lim}}}}

\newdimen\buffer@
\buffer@\fontdimen13 \tenex
\newdimen\buffer
\buffer\buffer@

\def\ResetBuffer{\fontdimen13 \tenex\buffer@\global\buffer\buffer@}
\def\shave#1{\mathop{\hbox{$\m@th\fontdimen13 \tenex\z@                     
 \displaystyle{#1}$}}\fontdimen13 \tenex\buffer}

\message{multilevel sub/superscripts,}
\Invalid@\\
\def\Let@{\relax\iffalse{\fi\let\\=\cr\iffalse}\fi}
\Invalid@\vspace
\def\vspace@{\def\vspace##1{\crcr\noalign{\vskip##1\relax}}}
\def\multilimits@{\bgroup\vspace@\Let@
 \baselineskip\fontdimen10 \scriptfont\tw@
 \advance\baselineskip\fontdimen12 \scriptfont\tw@
 \lineskip\thr@@\fontdimen8 \scriptfont\thr@@
 \lineskiplimit\lineskip
 \vbox\bgroup\ialign\bgroup\hfil$\m@th\scriptstyle{##}$\hfil\crcr}
\def\Sb{_\multilimits@}
\def\endSb{\crcr\egroup\egroup\egroup}
\def\Sp{^\multilimits@}

\def\spreadlines#1{\RIfMIfI@\onlydmatherr@\spreadlines\else
 \openup#1\relax\fi\else\onlydmatherr@\spreadlines\fi}
\def\Mathstrut@{\copy\Mathstrutbox@}
\newbox\Mathstrutbox@
\setbox\Mathstrutbox@\null
\setboxz@h{$\m@th($}
\ht\Mathstrutbox@\ht\z@
\dp\Mathstrutbox@\dp\z@
\message{matrices,}
\newdimen\spreadmlines@
\def\spreadmatrixlines#1{\RIfMIfI@
 \onlydmatherr@\spreadmatrixlines\else
 \spreadmlines@#1\relax\fi\else\onlydmatherr@\spreadmatrixlines\fi}
\def\matrix{\null\,\vcenter\bgroup\Let@\vspace@
 \normalbaselines\openup\spreadmlines@\ialign
 \bgroup\hfil$\m@th##$\hfil&&\quad\hfil$\m@th##$\hfil\crcr
 \Mathstrut@\crcr\noalign{\kern-\baselineskip}}
\def\endmatrix{\crcr\Mathstrut@\crcr\noalign{\kern-\baselineskip}\egroup
 \egroup\,}
\def\format{\crcr\egroup\iffalse{\fi\ifnum`}=0 \fi\format@}
\newtoks\hashtoks@
\hashtoks@{#}
\def\format@#1\\{\def\preamble@{#1}%
 \def\l{$\m@th\the\hashtoks@$\hfil}%
 \def\c{\hfil$\m@th\the\hashtoks@$\hfil}%
 \def\r{\hfil$\m@th\the\hashtoks@$}%
 \edef\preamble@@{\preamble@}\ifnum`{=0 \fi\iffalse}\fi
 \ialign\bgroup\span\preamble@@\crcr}
\def\smallmatrix{\null\,\vcenter\bgroup\vspace@\Let@
 \baselineskip9\ex@\lineskip\ex@
 \ialign\bgroup\hfil$\m@th\scriptstyle{##}$\hfil&&\thickspace\hfil
 $\m@th\scriptstyle{##}$\hfil\crcr}
\def\endsmallmatrix{\crcr\egroup\egroup\,}

\newmuskip\dotsspace@
\dotsspace@1.5mu
\def\strip@#1 {#1}
\def\spacehdots#1\for#2{\multispan{#2}\xleaders
 \hbox{$\m@th\mkern\strip@#1 \dotsspace@.\mkern\strip@#1 \dotsspace@$}\hfill}
\def\hdotsfor#1{\spacehdots\@ne\for{#1}}
\def\multispan@#1{\omit\mscount#1\unskip\loop\ifnum\mscount>\@ne\sp@n\repeat}
\def\spaceinnerhdots#1\for#2\after#3{\multispan@{\strip@#2 }#3\xleaders
 \hbox{$\m@th\mkern\strip@#1 \dotsspace@.\mkern\strip@#1 \dotsspace@$}\hfill}
\def\innerhdotsfor#1\after#2{\spaceinnerhdots\@ne\for#1\after{#2}}
\def\cases{\bgroup\spreadmlines@\jot\left\{\,\matrix\format\l&\quad\l\\}
\def\endcases{\endmatrix\right.\egroup}
\message{multiline displays,}
\newif\ifinany@
\newif\ifinalign@
\newif\ifingather@
\def\strut@{\copy\strutbox@}
\newbox\strutbox@
\setbox\strutbox@\hbox{\vrule height8\p@ depth3\p@ width\z@}
\def\topaligned{\null\,\vtop\aligned@}
\def\botaligned{\null\,\vbox\aligned@}
\def\aligned{\null\,\vcenter\aligned@}
\def\aligned@{\bgroup\vspace@\Let@
 \ifinany@\else\openup\jot\fi\ialign
 \bgroup\hfil\strut@$\m@th\displaystyle{##}$&
 $\m@th\displaystyle{{}##}$\hfil\crcr}
\def\endaligned{\crcr\egroup\egroup}

\def\alignedat#1{\null\,\vcenter\bgroup\doat@{#1}\vspace@\Let@
 \ifinany@\else\openup\jot\fi\ialign\bgroup\span\preamble@@\crcr}
\newcount\atcount@
\def\doat@#1{\toks@{\hfil\strut@$\m@th
 \displaystyle{\the\hashtoks@}$&$\m@th\displaystyle
 {{}\the\hashtoks@}$\hfil}
 \atcount@#1\relax\advance\atcount@\m@ne                                    
 \loop\ifnum\atcount@>\z@\toks@=\expandafter{\the\toks@&\hfil$\m@th
 \displaystyle{\the\hashtoks@}$&$\m@th
 \displaystyle{{}\the\hashtoks@}$\hfil}\advance
  \atcount@\m@ne\repeat                                                     
 \xdef\preamble@{\the\toks@}\xdef\preamble@@{\preamble@}}

\def\gathered{\null\,\vcenter\bgroup\vspace@\Let@
 \ifinany@\else\openup\jot\fi\ialign
 \bgroup\hfil\strut@$\m@th\displaystyle{##}$\hfil\crcr}
\def\endgathered{\crcr\egroup\egroup}
\newif\iftagsleft@
\def\TagsOnLeft{\global\tagsleft@true}
\def\TagsOnRight{\global\tagsleft@false}
\TagsOnLeft
\newif\ifmathtags@
\def\TagsAsMath{\global\mathtags@true}
\def\TagsAsText{\global\mathtags@false}
\TagsAsText
\def\tagform@#1{\hbox{\rm(\ignorespaces#1\unskip)}}
\def\thetag{\leavevmode\tagform@}
\def\tag#1$${\iftagsleft@\leqno\else\eqno\fi                                
 \maketag@#1\maketag@                                                       
 $$}                                                                        
\def\maketag@{\FN@\maketag@@}
\def\maketag@@{\ifx\next"\expandafter\maketag@@@\else\expandafter\maketag@@@@
 \fi}
\def\maketag@@@"#1"#2\maketag@{\hbox{\rm#1}}                                
\def\maketag@@@@#1\maketag@{\ifmathtags@\tagform@{$\m@th#1$}\else
 \tagform@{#1}\fi}
\interdisplaylinepenalty\@M
\def\allowdisplaybreaks{\RIfMIfI@
 \onlydmatherr@\allowdisplaybreaks\else
 \interdisplaylinepenalty\z@\fi\else\onlydmatherr@\allowdisplaybreaks\fi}
\Invalid@\allowdisplaybreak
\Invalid@\displaybreak
\Invalid@\intertext
\def\allowdisplaybreak@{\def\allowdisplaybreak{\crcr\noalign{\allowbreak}}}
\def\displaybreak@{\def\displaybreak{\crcr\noalign{\break}}}
\def\intertext@{\def\intertext##1{\crcr\noalign{%
 \penalty\postdisplaypenalty \vskip\belowdisplayskip
 \vbox{\normalbaselines\noindent##1}%
 \penalty\predisplaypenalty \vskip\abovedisplayskip}}}
\newskip\centering@
\centering@\z@ plus\@m\p@
\def\align{\relax\ifingather@\DN@{\csname align (in
  \string\gather)\endcsname}\else
 \ifmmode\ifinner\DN@{\onlydmatherr@\align}\else
  \let\next@\align@\fi
 \else\DN@{\onlydmatherr@\align}\fi\fi\next@}
\newhelp\andhelp@
{An extra & here is so disastrous that you should probably exit^^J
and fix things up.}
\newif\iftag@
\newcount\and@
\def\align@{\inalign@true\inany@true
 \vspace@\allowdisplaybreak@\displaybreak@\intertext@
 \def\tag{\global\tag@true\ifnum\and@=\z@\DN@{&&}\else
          \DN@{&}\fi\next@}%
 \iftagsleft@\DN@{\csname align \endcsname}\else
  \DN@{\csname align \space\endcsname}\fi\next@}
\def\Tag@{\iftag@\else\errhelp\andhelp@\err@{Extra & on this line}\fi}
\newdimen\lwidth@
\newdimen\rwidth@
\newdimen\maxlwidth@
\newdimen\maxrwidth@
\newdimen\totwidth@
\def\measure@#1\endalign{\lwidth@\z@\rwidth@\z@\maxlwidth@\z@\maxrwidth@\z@
 \global\and@\z@                                                            
 \setbox@ne\vbox                                                            
  {\everycr{\noalign{\global\tag@false\global\and@\z@}}\Let@                
  \halign{\setboxz@h{$\m@th\displaystyle{\@lign##}$}
   \global\lwidth@\wdz@                                                     
   \ifdim\lwidth@>\maxlwidth@\global\maxlwidth@\lwidth@\fi                  
   \global\advance\and@\@ne                                                 
   &\setboxz@h{$\m@th\displaystyle{{}\@lign##}$}\global\rwidth@\wdz@        
   \ifdim\rwidth@>\maxrwidth@\global\maxrwidth@\rwidth@\fi                  
   \global\advance\and@\@ne                                                
   &\Tag@
   \eat@{##}\crcr#1\crcr}}
 \totwidth@\maxlwidth@\advance\totwidth@\maxrwidth@}                       
\def\displ@y@{\global\dt@ptrue\openup\jot
 \everycr{\noalign{\global\tag@false\global\and@\z@\ifdt@p\global\dt@pfalse
 \vskip-\lineskiplimit\vskip\normallineskiplimit\else
 \penalty\interdisplaylinepenalty\fi}}}
\def\black@#1{\noalign{\ifdim#1>\displaywidth
 \dimen@\prevdepth\nointerlineskip                                          
 \vskip-\ht\strutbox@\vskip-\dp\strutbox@                                   
 \vbox{\noindent\hbox to#1{\strut@\hfill}}
 \prevdepth\dimen@                                                          
 \fi}}
\expandafter\def\csname align \space\endcsname#1\endalign
 {\measure@#1\endalign\global\and@\z@                                       
 \ifingather@\everycr{\noalign{\global\and@\z@}}\else\displ@y@\fi           
 \Let@\tabskip\centering@                                                   
 \halign to\displaywidth
  {\hfil\strut@\setboxz@h{$\m@th\displaystyle{\@lign##}$}
  \global\lwidth@\wdz@\boxz@\global\advance\and@\@ne                        
  \tabskip\z@skip                                                           
  &\setboxz@h{$\m@th\displaystyle{{}\@lign##}$}
  \global\rwidth@\wdz@\boxz@\hfill\global\advance\and@\@ne                  
  \tabskip\centering@                                                       
  &\setboxz@h{\@lign\strut@\maketag@##\maketag@}
  \dimen@\displaywidth\advance\dimen@-\totwidth@
  \divide\dimen@\tw@\advance\dimen@\maxrwidth@\advance\dimen@-\rwidth@     
  \ifdim\dimen@<\tw@\wdz@\llap{\vtop{\normalbaselines\null\boxz@}}
  \else\llap{\boxz@}\fi                                                    
  \tabskip\z@skip                                                          
  \crcr#1\crcr                                                             
  \black@\totwidth@}}                                                      
\newdimen\lineht@
\expandafter\def\csname align \endcsname#1\endalign{\measure@#1\endalign
 \global\and@\z@
 \ifdim\totwidth@>\displaywidth\let\displaywidth@\totwidth@\else
  \let\displaywidth@\displaywidth\fi                                        
 \ifingather@\everycr{\noalign{\global\and@\z@}}\else\displ@y@\fi
 \Let@\tabskip\centering@\halign to\displaywidth
  {\hfil\strut@\setboxz@h{$\m@th\displaystyle{\@lign##}$}%
  \global\lwidth@\wdz@\global\lineht@\ht\z@                                 
  \boxz@\global\advance\and@\@ne
  \tabskip\z@skip&\setboxz@h{$\m@th\displaystyle{{}\@lign##}$}%
  \global\rwidth@\wdz@\ifdim\ht\z@>\lineht@\global\lineht@\ht\z@\fi         
  \boxz@\hfil\global\advance\and@\@ne
  \tabskip\centering@&\kern-\displaywidth@                                  
  \setboxz@h{\@lign\strut@\maketag@##\maketag@}%
  \dimen@\displaywidth\advance\dimen@-\totwidth@
  \divide\dimen@\tw@\advance\dimen@\maxlwidth@\advance\dimen@-\lwidth@
  \ifdim\dimen@<\tw@\wdz@
   \rlap{\vbox{\normalbaselines\boxz@\vbox to\lineht@{}}}\else
   \rlap{\boxz@}\fi
  \tabskip\displaywidth@\crcr#1\crcr\black@\totwidth@}}
\expandafter\def\csname align (in \string\gather)\endcsname
  #1\endalign{\vcenter{\align@#1\endalign}}
\Invalid@\endalign
\newif\ifxat@
\def\alignat{\RIfMIfI@\DN@{\onlydmatherr@\alignat}\else
 \DN@{\csname alignat \endcsname}\fi\else
 \DN@{\onlydmatherr@\alignat}\fi\next@}
\newif\ifmeasuring@
\newbox\savealignat@
\expandafter\def\csname alignat \endcsname#1#2\endalignat                   
 {\inany@true\xat@false
 \def\tag{\global\tag@true\count@#1\relax\multiply\count@\tw@
  \xdef\tag@{}\loop\ifnum\count@>\and@\xdef\tag@{&\tag@}\advance\count@\m@ne
  \repeat\tag@}%
 \vspace@\allowdisplaybreak@\displaybreak@\intertext@
 \displ@y@\measuring@true                                                   
 \setbox\savealignat@\hbox{$\m@th\displaystyle\Let@
  \attag@{#1}
  \vbox{\halign{\span\preamble@@\crcr#2\crcr}}$}%
 \measuring@false                                                           
 \Let@\attag@{#1}
 \tabskip\centering@\halign to\displaywidth
  {\span\preamble@@\crcr#2\crcr                                             
  \black@{\wd\savealignat@}}}                                               
\Invalid@\endalignat
\def\xalignat{\RIfMIfI@
 \DN@{\onlydmatherr@\xalignat}\else
 \DN@{\csname xalignat \endcsname}\fi\else
 \DN@{\onlydmatherr@\xalignat}\fi\next@}
\expandafter\def\csname xalignat \endcsname#1#2\endxalignat
 {\inany@true\xat@true
 \def\tag{\global\tag@true\def\tag@{}\count@#1\relax\multiply\count@\tw@
  \loop\ifnum\count@>\and@\xdef\tag@{&\tag@}\advance\count@\m@ne\repeat\tag@}%
 \vspace@\allowdisplaybreak@\displaybreak@\intertext@
 \displ@y@\measuring@true\setbox\savealignat@\hbox{$\m@th\displaystyle\Let@
 \attag@{#1}\vbox{\halign{\span\preamble@@\crcr#2\crcr}}$}%
 \measuring@false\Let@
 \attag@{#1}\tabskip\centering@\halign to\displaywidth
 {\span\preamble@@\crcr#2\crcr\black@{\wd\savealignat@}}}
\def\attag@#1{\let\Maketag@\maketag@\let\TAG@\Tag@                          
 \let\Tag@=0\let\maketag@=0
 \ifmeasuring@\def\llap@##1{\setboxz@h{##1}\hbox to\tw@\wdz@{}}%
  \def\rlap@##1{\setboxz@h{##1}\hbox to\tw@\wdz@{}}\else
  \let\llap@\llap\let\rlap@\rlap\fi                                         
 \toks@{\hfil\strut@$\m@th\displaystyle{\@lign\the\hashtoks@}$\tabskip\z@skip
  \global\advance\and@\@ne&$\m@th\displaystyle{{}\@lign\the\hashtoks@}$\hfil
  \ifxat@\tabskip\centering@\fi\global\advance\and@\@ne}
 \iftagsleft@
  \toks@@{\tabskip\centering@&\Tag@\kern-\displaywidth
   \rlap@{\@lign\maketag@\the\hashtoks@\maketag@}%
   \global\advance\and@\@ne\tabskip\displaywidth}\else
  \toks@@{\tabskip\centering@&\Tag@\llap@{\@lign\maketag@
   \the\hashtoks@\maketag@}\global\advance\and@\@ne\tabskip\z@skip}\fi      
 \atcount@#1\relax\advance\atcount@\m@ne
 \loop\ifnum\atcount@>\z@
 \toks@=\expandafter{\the\toks@&\hfil$\m@th\displaystyle{\@lign
  \the\hashtoks@}$\global\advance\and@\@ne
  \tabskip\z@skip&$\m@th\displaystyle{{}\@lign\the\hashtoks@}$\hfil\ifxat@
  \tabskip\centering@\fi\global\advance\and@\@ne}\advance\atcount@\m@ne
 \repeat                                                                    
 \xdef\preamble@{\the\toks@\the\toks@@}
 \xdef\preamble@@{\preamble@}
 \let\maketag@\Maketag@\let\Tag@\TAG@}                                      
\Invalid@\endxalignat
\def\xxalignat{\RIfMIfI@
 \DN@{\onlydmatherr@\xxalignat}\else\DN@{\csname xxalignat
  \endcsname}\fi\else
 \DN@{\onlydmatherr@\xxalignat}\fi\next@}
\expandafter\def\csname xxalignat \endcsname#1#2\endxxalignat{\inany@true
 \vspace@\allowdisplaybreak@\displaybreak@\intertext@
 \displ@y\setbox\savealignat@\hbox{$\m@th\displaystyle\Let@
 \xxattag@{#1}\vbox{\halign{\span\preamble@@\crcr#2\crcr}}$}%
 \Let@\xxattag@{#1}\tabskip\z@skip\halign to\displaywidth
 {\span\preamble@@\crcr#2\crcr\black@{\wd\savealignat@}}}
\def\xxattag@#1{\toks@{\tabskip\z@skip\hfil\strut@
 $\m@th\displaystyle{\the\hashtoks@}$&%
 $\m@th\displaystyle{{}\the\hashtoks@}$\hfil\tabskip\centering@&}%
 \atcount@#1\relax\advance\atcount@\m@ne\loop\ifnum\atcount@>\z@
 \toks@=\expandafter{\the\toks@&\hfil$\m@th\displaystyle{\the\hashtoks@}$%
  \tabskip\z@skip&$\m@th\displaystyle{{}\the\hashtoks@}$\hfil
  \tabskip\centering@}\advance\atcount@\m@ne\repeat
 \xdef\preamble@{\the\toks@\tabskip\z@skip}\xdef\preamble@@{\preamble@}}
\Invalid@\endxxalignat
\newdimen\gwidth@
\newdimen\gmaxwidth@
\def\gmeasure@#1\endgather{\gwidth@\z@\gmaxwidth@\z@\setbox@ne\vbox{\Let@
 \halign{\setboxz@h{$\m@th\displaystyle{##}$}\global\gwidth@\wdz@
 \ifdim\gwidth@>\gmaxwidth@\global\gmaxwidth@\gwidth@\fi
 &\eat@{##}\crcr#1\crcr}}}
\def\gather{\RIfMIfI@\DN@{\onlydmatherr@\gather}\else
 \ingather@true\inany@true\def\tag{&}%
 \vspace@\allowdisplaybreak@\displaybreak@\intertext@
 \displ@y\Let@
 \iftagsleft@\DN@{\csname gather \endcsname}\else
  \DN@{\csname gather \space\endcsname}\fi\fi
 \else\DN@{\onlydmatherr@\gather}\fi\next@}
\expandafter\def\csname gather \space\endcsname#1\endgather
 {\gmeasure@#1\endgather\tabskip\centering@
 \halign to\displaywidth{\hfil\strut@\setboxz@h{$\m@th\displaystyle{##}$}%
 \global\gwidth@\wdz@\boxz@\hfil&
 \setboxz@h{\strut@{\maketag@##\maketag@}}%
 \dimen@\displaywidth\advance\dimen@-\gwidth@
 \ifdim\dimen@>\tw@\wdz@\llap{\boxz@}\else
 \llap{\vtop{\normalbaselines\null\boxz@}}\fi
 \tabskip\z@skip\crcr#1\crcr\black@\gmaxwidth@}}
\newdimen\glineht@
\expandafter\def\csname gather \endcsname#1\endgather{\gmeasure@#1\endgather
 \ifdim\gmaxwidth@>\displaywidth\let\gdisplaywidth@\gmaxwidth@\else
 \let\gdisplaywidth@\displaywidth\fi\tabskip\centering@\halign to\displaywidth
 {\hfil\strut@\setboxz@h{$\m@th\displaystyle{##}$}%
 \global\gwidth@\wdz@\global\glineht@\ht\z@\boxz@\hfil&\kern-\gdisplaywidth@
 \setboxz@h{\strut@{\maketag@##\maketag@}}%
 \dimen@\displaywidth\advance\dimen@-\gwidth@
 \ifdim\dimen@>\tw@\wdz@\rlap{\boxz@}\else
 \rlap{\vbox{\normalbaselines\boxz@\vbox to\glineht@{}}}\fi
 \tabskip\gdisplaywidth@\crcr#1\crcr\black@\gmaxwidth@}}
\newif\ifctagsplit@
\def\CenteredTagsOnSplits{\global\ctagsplit@true}
\def\TopOrBottomTagsOnSplits{\global\ctagsplit@false}
\TopOrBottomTagsOnSplits
\def\split{\relax\ifinany@\let\next@\insplit@\else
 \ifmmode\ifinner\def\next@{\onlydmatherr@\split}\else
 \let\next@\outsplit@\fi\else
 \def\next@{\onlydmatherr@\split}\fi\fi\next@}
\def\insplit@{\global\setbox\z@\vbox\bgroup\vspace@\Let@\ialign\bgroup
 \hfil\strut@$\m@th\displaystyle{##}$&$\m@th\displaystyle{{}##}$\hfill\crcr}
\def\endsplit{\crcr\egroup\egroup\iftagsleft@\expandafter\lendsplit@\else
 \expandafter\rendsplit@\fi}
\def\rendsplit@{\global\setbox9 \vbox
 {\unvcopy\z@\global\setbox8 \lastbox\unskip}
 \setbox@ne\hbox{\unhcopy8 \unskip\global\setbox\tw@\lastbox
 \unskip\global\setbox\thr@@\lastbox}
 \global\setbox7 \hbox{\unhbox\tw@\unskip}
 \ifinalign@\ifctagsplit@                                                   
  \gdef\split@{\hbox to\wd\thr@@{}&
   \vcenter{\vbox{\moveleft\wd\thr@@\boxz@}}}
 \else\gdef\split@{&\vbox{\moveleft\wd\thr@@\box9}\crcr
  \box\thr@@&\box7}\fi                                                      
 \else                                                                      
  \ifctagsplit@\gdef\split@{\vcenter{\boxz@}}\else
  \gdef\split@{\box9\crcr\hbox{\box\thr@@\box7}}\fi
 \fi
 \split@}                                                                   
\def\lendsplit@{\global\setbox9\vtop{\unvcopy\z@}
 \setbox@ne\vbox{\unvcopy\z@\global\setbox8\lastbox}
 \setbox@ne\hbox{\unhcopy8\unskip\setbox\tw@\lastbox
  \unskip\global\setbox\thr@@\lastbox}
 \ifinalign@\ifctagsplit@                                                   
  \gdef\split@{\hbox to\wd\thr@@{}&
  \vcenter{\vbox{\moveleft\wd\thr@@\box9}}}
  \else                                                                     
  \gdef\split@{\hbox to\wd\thr@@{}&\vbox{\moveleft\wd\thr@@\box9}}\fi
 \else
  \ifctagsplit@\gdef\split@{\vcenter{\box9}}\else
  \gdef\split@{\box9}\fi
 \fi\split@}
\def\outsplit@#1$${\align\insplit@#1\endalign$$}
\newdimen\multlinegap@
\multlinegap@1em
\newdimen\multlinetaggap@
\multlinetaggap@1em
\def\MultlineGap#1{\global\multlinegap@#1\relax}
\def\multlinegap#1{\RIfMIfI@\onlydmatherr@\multlinegap\else
 \multlinegap@#1\relax\fi\else\onlydmatherr@\multlinegap\fi}
\def\nomultlinegap{\multlinegap{\z@}}
\def\multline{\RIfMIfI@
 \DN@{\onlydmatherr@\multline}\else
 \DN@{\multline@}\fi\else
 \DN@{\onlydmatherr@\multline}\fi\next@}
\newif\iftagin@
\def\tagin@#1{\tagin@false\in@\tag{#1}\ifin@\tagin@true\fi}
\def\multline@#1$${\inany@true\vspace@\allowdisplaybreak@\displaybreak@
 \tagin@{#1}\iftagsleft@\DN@{\multline@l#1$$}\else
 \DN@{\multline@r#1$$}\fi\next@}
\newdimen\mwidth@
\def\rmmeasure@#1\endmultline{%
 \def\shoveleft##1{##1}\def\shoveright##1{##1}
 \setbox@ne\vbox{\Let@\halign{\setboxz@h
  {$\m@th\@lign\displaystyle{}##$}\global\mwidth@\wdz@
  \crcr#1\crcr}}}
\newdimen\mlineht@
\newif\ifzerocr@
\newif\ifonecr@
\def\lmmeasure@#1\endmultline{\global\zerocr@true\global\onecr@false
 \everycr{\noalign{\ifonecr@\global\onecr@false\fi
  \ifzerocr@\global\zerocr@false\global\onecr@true\fi}}
  \def\shoveleft##1{##1}\def\shoveright##1{##1}%
 \setbox@ne\vbox{\Let@\halign{\setboxz@h
  {$\m@th\@lign\displaystyle{}##$}\ifonecr@\global\mwidth@\wdz@
  \global\mlineht@\ht\z@\fi\crcr#1\crcr}}}
\newbox\mtagbox@
\newdimen\ltwidth@
\newdimen\rtwidth@
\def\multline@l#1$${\iftagin@\DN@{\lmultline@@#1$$}\else
 \DN@{\setbox\mtagbox@\null\ltwidth@\z@\rtwidth@\z@
  \lmultline@@@#1$$}\fi\next@}
\def\lmultline@@#1\endmultline\tag#2$${%
 \setbox\mtagbox@\hbox{\maketag@#2\maketag@}
 \lmmeasure@#1\endmultline\dimen@\mwidth@\advance\dimen@\wd\mtagbox@
 \advance\dimen@\multlinetaggap@                                            
 \ifdim\dimen@>\displaywidth\ltwidth@\z@\else\ltwidth@\wd\mtagbox@\fi       
 \lmultline@@@#1\endmultline$$}
\def\lmultline@@@{\displ@y
 \def\shoveright##1{##1\hfilneg\hskip\multlinegap@}%
 \def\shoveleft##1{\setboxz@h{$\m@th\displaystyle{}##1$}%
  \setbox@ne\hbox{$\m@th\displaystyle##1$}%
  \hfilneg
  \iftagin@
   \ifdim\ltwidth@>\z@\hskip\ltwidth@\hskip\multlinetaggap@\fi
  \else\hskip\multlinegap@\fi\hskip.5\wd@ne\hskip-.5\wdz@##1}
  \halign\bgroup\Let@\hbox to\displaywidth
   {\strut@$\m@th\displaystyle\hfil{}##\hfil$}\crcr
   \hfilneg                                                                 
   \iftagin@                                                                
    \ifdim\ltwidth@>\z@                                                     
     \box\mtagbox@\hskip\multlinetaggap@                                    
    \else
     \rlap{\vbox{\normalbaselines\hbox{\strut@\box\mtagbox@}%
     \vbox to\mlineht@{}}}\fi                                               
   \else\hskip\multlinegap@\fi}                                             
\def\multline@r#1$${\iftagin@\DN@{\rmultline@@#1$$}\else
 \DN@{\setbox\mtagbox@\null\ltwidth@\z@\rtwidth@\z@
  \rmultline@@@#1$$}\fi\next@}
\def\rmultline@@#1\endmultline\tag#2$${\ltwidth@\z@
 \setbox\mtagbox@\hbox{\maketag@#2\maketag@}%
 \rmmeasure@#1\endmultline\dimen@\mwidth@\advance\dimen@\wd\mtagbox@
 \advance\dimen@\multlinetaggap@
 \ifdim\dimen@>\displaywidth\rtwidth@\z@\else\rtwidth@\wd\mtagbox@\fi
 \rmultline@@@#1\endmultline$$}
\def\rmultline@@@{\displ@y
 \def\shoveright##1{##1\hfilneg\iftagin@\ifdim\rtwidth@>\z@
  \hskip\rtwidth@\hskip\multlinetaggap@\fi\else\hskip\multlinegap@\fi}%
 \def\shoveleft##1{\setboxz@h{$\m@th\displaystyle{}##1$}%
  \setbox@ne\hbox{$\m@th\displaystyle##1$}%
  \hfilneg\hskip\multlinegap@\hskip.5\wd@ne\hskip-.5\wdz@##1}%
 \halign\bgroup\Let@\hbox to\displaywidth
  {\strut@$\m@th\displaystyle\hfil{}##\hfil$}\crcr
 \hfilneg\hskip\multlinegap@}
\def\endmultline{\iftagsleft@\expandafter\lendmultline@\else
 \expandafter\rendmultline@\fi}
\def\lendmultline@{\hfilneg\hskip\multlinegap@\crcr\egroup}
\def\rendmultline@{\iftagin@                                                
 \ifdim\rtwidth@>\z@                                                        
  \hskip\multlinetaggap@\box\mtagbox@                                       
 \else\llap{\vtop{\normalbaselines\null\hbox{\strut@\box\mtagbox@}}}\fi     
 \else\hskip\multlinegap@\fi                                                
 \hfilneg\crcr\egroup}
\def\bmod{\mskip-\medmuskip\mkern5mu\mathbin{\fam\z@ mod}\penalty900
 \mkern5mu\mskip-\medmuskip}
\def\pmod#1{\allowbreak\ifinner\mkern8mu\else\mkern18mu\fi
 ({\fam\z@ mod}\,\,#1)}
\def\pod#1{\allowbreak\ifinner\mkern8mu\else\mkern18mu\fi(#1)}
\def\mod#1{\allowbreak\ifinner\mkern12mu\else\mkern18mu\fi{\fam\z@ mod}\,\,#1}
\message{continued fractions,}
\newcount\cfraccount@
\def\cfrac{\bgroup\bgroup\advance\cfraccount@\@ne\strut
 \iffalse{\fi\def\\{\over\displaystyle}\iffalse}\fi}
\def\lcfrac{\bgroup\bgroup\advance\cfraccount@\@ne\strut
 \iffalse{\fi\def\\{\hfill\over\displaystyle}\iffalse}\fi}
\def\rcfrac{\bgroup\bgroup\advance\cfraccount@\@ne\strut\hfill
 \iffalse{\fi\def\\{\over\displaystyle}\iffalse}\fi}
\def\gloop@#1\repeat{\gdef\body{#1}\iterate}
\def\endcfrac{\gloop@\ifnum\cfraccount@>\z@\global\advance\cfraccount@\m@ne
 \egroup\hskip-\nulldelimiterspace\egroup\repeat}
\message{compound symbols,}
\def\binrel@#1{\setboxz@h{\thinmuskip0mu
  \medmuskip\m@ne mu\thickmuskip\@ne mu$#1\m@th$}%
 \setbox@ne\hbox{\thinmuskip0mu\medmuskip\m@ne mu\thickmuskip
  \@ne mu${}#1{}\m@th$}%
 \setbox\tw@\hbox{\hskip\wd@ne\hskip-\wdz@}}
\def\overset#1\to#2{\binrel@{#2}\ifdim\wd\tw@<\z@
 \mathbin{\mathop{\kern\z@#2}\limits^{#1}}\else\ifdim\wd\tw@>\z@
 \mathrel{\mathop{\kern\z@#2}\limits^{#1}}\else
 {\mathop{\kern\z@#2}\limits^{#1}}{}\fi\fi}
\def\underset#1\to#2{\binrel@{#2}\ifdim\wd\tw@<\z@
 \mathbin{\mathop{\kern\z@#2}\limits_{#1}}\else\ifdim\wd\tw@>\z@
 \mathrel{\mathop{\kern\z@#2}\limits_{#1}}\else
 {\mathop{\kern\z@#2}\limits_{#1}}{}\fi\fi}
\def\oversetbrace#1\to#2{\overbrace{#2}^{#1}}
\def\undersetbrace#1\to#2{\underbrace{#2}_{#1}}
\def\sideset#1\and#2\to#3{%
 \setbox@ne\hbox{$\dsize{\vphantom{#3}}#1{#3}\m@th$}%
 \setbox\tw@\hbox{$\dsize{#3}#2\m@th$}%
 \hskip\wd@ne\hskip-\wd\tw@\mathop{\hskip\wd\tw@\hskip-\wd@ne
  {\vphantom{#3}}#1{#3}#2}}
\def\rightarrowfill@#1{\setboxz@h{$#1-\m@th$}\ht\z@\z@
  $#1\m@th\copy\z@\mkern-6mu\cleaders
  \hbox{$#1\mkern-2mu\box\z@\mkern-2mu$}\hfill
  \mkern-6mu\mathord\rightarrow$}
\def\leftarrowfill@#1{\setboxz@h{$#1-\m@th$}\ht\z@\z@
  $#1\m@th\mathord\leftarrow\mkern-6mu\cleaders
  \hbox{$#1\mkern-2mu\copy\z@\mkern-2mu$}\hfill
  \mkern-6mu\box\z@$}
\def\leftrightarrowfill@#1{\setboxz@h{$#1-\m@th$}\ht\z@\z@
  $#1\m@th\mathord\leftarrow\mkern-6mu\cleaders
  \hbox{$#1\mkern-2mu\box\z@\mkern-2mu$}\hfill
  \mkern-6mu\mathord\rightarrow$}
\def\overrightarrow{\mathpalette\overrightarrow@}
\def\overrightarrow@#1#2{\vbox{\ialign{##\crcr\rightarrowfill@#1\crcr
 \noalign{\kern-\ex@\nointerlineskip}$\m@th\hfil#1#2\hfil$\crcr}}}

\def\overleftarrow{\mathpalette\overleftarrow@}
\def\overleftarrow@#1#2{\vbox{\ialign{##\crcr\leftarrowfill@#1\crcr
 \noalign{\kern-\ex@\nointerlineskip}$\m@th\hfil#1#2\hfil$\crcr}}}
\def\overleftrightarrow{\mathpalette\overleftrightarrow@}
\def\overleftrightarrow@#1#2{\vbox{\ialign{##\crcr\leftrightarrowfill@#1\crcr
 \noalign{\kern-\ex@\nointerlineskip}$\m@th\hfil#1#2\hfil$\crcr}}}
\def\underrightarrow{\mathpalette\underrightarrow@}
\def\underrightarrow@#1#2{\vtop{\ialign{##\crcr$\m@th\hfil#1#2\hfil$\crcr
 \noalign{\nointerlineskip}\rightarrowfill@#1\crcr}}}

\def\underleftarrow{\mathpalette\underleftarrow@}
\def\underleftarrow@#1#2{\vtop{\ialign{##\crcr$\m@th\hfil#1#2\hfil$\crcr
 \noalign{\nointerlineskip}\leftarrowfill@#1\crcr}}}
\def\underleftrightarrow{\mathpalette\underleftrightarrow@}
\def\underleftrightarrow@#1#2{\vtop{\ialign{##\crcr$\m@th\hfil#1#2\hfil$\crcr
 \noalign{\nointerlineskip}\leftrightarrowfill@#1\crcr}}}
\message{various kinds of dots,}
\let\DOTSI\relax
\let\DOTSB\relax

\newif\ifmath@
{\uccode`7=`\\ \uccode`8=`m \uccode`9=`a \uccode`0=`t \uccode`!=`h
 \uppercase{\gdef\math@#1#2#3#4#5#6\math@{\global\math@false\ifx 7#1\ifx 8#2%
 \ifx 9#3\ifx 0#4\ifx !#5\xdef\meaning@{#6}\global\math@true\fi\fi\fi\fi\fi}}}
\newif\ifmathch@
{\uccode`7=`c \uccode`8=`h \uccode`9=`\"
 \uppercase{\gdef\mathch@#1#2#3#4#5#6\mathch@{\global\mathch@false
  \ifx 7#1\ifx 8#2\ifx 9#5\global\mathch@true\xdef\meaning@{9#6}\fi\fi\fi}}}
\newcount\classnum@
\def\getmathch@#1.#2\getmathch@{\classnum@#1 \divide\classnum@4096
 \ifcase\number\classnum@\or\or\gdef\thedots@{\dotsb@}\or
 \gdef\thedots@{\dotsb@}\fi}
\newif\ifmathbin@
{\uccode`4=`b \uccode`5=`i \uccode`6=`n
 \uppercase{\gdef\mathbin@#1#2#3{\relaxnext@
  \DNii@##1\mathbin@{\ifx\space@\next\global\mathbin@true\fi}%
 \global\mathbin@false\DN@##1\mathbin@{}%
 \ifx 4#1\ifx 5#2\ifx 6#3\DN@{\FN@\nextii@}\fi\fi\fi\next@}}}
\newif\ifmathrel@
{\uccode`4=`r \uccode`5=`e \uccode`6=`l
 \uppercase{\gdef\mathrel@#1#2#3{\relaxnext@
  \DNii@##1\mathrel@{\ifx\space@\next\global\mathrel@true\fi}%
 \global\mathrel@false\DN@##1\mathrel@{}%
 \ifx 4#1\ifx 5#2\ifx 6#3\DN@{\FN@\nextii@}\fi\fi\fi\next@}}}
\newif\ifmacro@
{\uccode`5=`m \uccode`6=`a \uccode`7=`c
 \uppercase{\gdef\macro@#1#2#3#4\macro@{\global\macro@false
  \ifx 5#1\ifx 6#2\ifx 7#3\global\macro@true
  \xdef\meaning@{\macro@@#4\macro@@}\fi\fi\fi}}}
\def\macro@@#1->#2\macro@@{#2}
\newif\ifDOTS@
\newcount\DOTSCASE@
{\uccode`6=`\\ \uccode`7=`D \uccode`8=`O \uccode`9=`T \uccode`0=`S
 \uppercase{\gdef\DOTS@#1#2#3#4#5{\global\DOTS@false\DN@##1\DOTS@{}%
  \ifx 6#1\ifx 7#2\ifx 8#3\ifx 9#4\ifx 0#5\let\next@\DOTS@@\fi\fi\fi\fi\fi
  \next@}}}
{\uccode`3=`B \uccode`4=`I \uccode`5=`X
 \uppercase{\gdef\DOTS@@#1{\relaxnext@
  \DNii@##1\DOTS@{\ifx\space@\next\global\DOTS@true\fi}%
  \DN@{\FN@\nextii@}%
  \ifx 3#1\global\DOTSCASE@\z@\else
  \ifx 4#1\global\DOTSCASE@\@ne\else
  \ifx 5#1\global\DOTSCASE@\tw@\else\DN@##1\DOTS@{}%
  \fi\fi\fi\next@}}}
\newif\ifnot@
{\uccode`5=`\\ \uccode`6=`n \uccode`7=`o \uccode`8=`t
 \uppercase{\gdef\not@#1#2#3#4{\relaxnext@
  \DNii@##1\not@{\ifx\space@\next\global\not@true\fi}%
 \global\not@false\DN@##1\not@{}%
 \ifx 5#1\ifx 6#2\ifx 7#3\ifx 8#4\DN@{\FN@\nextii@}\fi\fi\fi
 \fi\next@}}}
\newif\ifkeybin@
\def\keybin@{\keybin@true
 \ifx\next+\else\ifx\next=\else\ifx\next<\else\ifx\next>\else\ifx\next-\else
 \ifx\next*\else\ifx\next:\else\keybin@false\fi\fi\fi\fi\fi\fi\fi}
\def\dots{\RIfM@\expandafter\mdots@\else\expandafter\tdots@\fi}
\def\tdots@{\unskip\relaxnext@
 \DN@{$\m@th\mathinner{\ldotp\ldotp\ldotp}\,
   \ifx\next,\,$\else\ifx\next.\,$\else\ifx\next;\,$\else\ifx\next:\,$\else
   \ifx\next?\,$\else\ifx\next!\,$\else$ \fi\fi\fi\fi\fi\fi}%
 \ \FN@\next@}
\def\mdots@{\FN@\mdots@@}
\def\mdots@@{\gdef\thedots@{\dotso@}
 \ifx\next\boldkey\gdef\thedots@\boldkey{\boldkeydots@}\else                
 \ifx\next\boldsymbol\gdef\thedots@\boldsymbol{\boldsymboldots@}\else       
 \ifx,\next\gdef\thedots@{\dotsc}
 \else\ifx\not\next\gdef\thedots@{\dotsb@}
 \else\keybin@
 \ifkeybin@\gdef\thedots@{\dotsb@}
 \else\xdef\meaning@{\meaning\next..........}\xdef\meaning@@{\meaning@}
  \expandafter\math@\meaning@\math@
  \ifmath@
   \expandafter\mathch@\meaning@\mathch@
   \ifmathch@\expandafter\getmathch@\meaning@\getmathch@\fi                 
  \else\expandafter\macro@\meaning@@\macro@                                 
  \ifmacro@                                                                
   \expandafter\not@\meaning@\not@\ifnot@\gdef\thedots@{\dotsb@}
  \else\expandafter\DOTS@\meaning@\DOTS@
  \ifDOTS@
   \ifcase\number\DOTSCASE@\gdef\thedots@{\dotsb@}%
    \or\gdef\thedots@{\dotsi}\else\fi                                      
  \else\expandafter\math@\meaning@\math@                                   
  \ifmath@\expandafter\mathbin@\meaning@\mathbin@
  \ifmathbin@\gdef\thedots@{\dotsb@}
  \else\expandafter\mathrel@\meaning@\mathrel@
  \ifmathrel@\gdef\thedots@{\dotsb@}
  \fi\fi\fi\fi\fi\fi\fi\fi\fi\fi\fi\fi
 \thedots@}
\def\plainldots@{\mathinner{\ldotp\ldotp\ldotp}}
\def\plaincdots@{\mathinner{\cdotp\cdotp\cdotp}}
\def\dotsi{\!\plaincdots@}
\let\dotsb@\plaincdots@
\newif\ifextra@
\newif\ifrightdelim@
\def\rightdelim@{\global\rightdelim@true                                    
 \ifx\next)\else                                                            
 \ifx\next]\else
 \ifx\next\rbrack\else
 \ifx\next\}\else
 \ifx\next\rbrace\else
 \ifx\next\rangle\else
 \ifx\next\rceil\else
 \ifx\next\rfloor\else
 \ifx\next\rgroup\else
 \ifx\next\rmoustache\else
 \ifx\next\right\else
 \ifx\next\bigr\else
 \ifx\next\biggr\else
 \ifx\next\Bigr\else                                                        
 \ifx\next\Biggr\else\global\rightdelim@false
 \fi\fi\fi\fi\fi\fi\fi\fi\fi\fi\fi\fi\fi\fi\fi}
\def\extra@{%
 \global\extra@false\rightdelim@\ifrightdelim@\global\extra@true            
 \else\ifx\next$\global\extra@true                                          
 \else\xdef\meaning@{\meaning\next..........}
 \expandafter\macro@\meaning@\macro@\ifmacro@                               
 \expandafter\DOTS@\meaning@\DOTS@
 \ifDOTS@
 \ifnum\DOTSCASE@=\tw@\global\extra@true                                    
 \fi\fi\fi\fi\fi}
\newif\ifbold@
\def\dotso@{\relaxnext@
 \ifbold@
  \let\next\delayed@
  \DNii@{\extra@\plainldots@\ifextra@\,\fi}%
 \else
  \DNii@{\DN@{\extra@\plainldots@\ifextra@\,\fi}\FN@\next@}%
 \fi
 \nextii@}
\def\extrap@#1{%
 \ifx\next,\DN@{#1\,}\else
 \ifx\next;\DN@{#1\,}\else
 \ifx\next.\DN@{#1\,}\else\extra@
 \ifextra@\DN@{#1\,}\else
 \let\next@#1\fi\fi\fi\fi\next@}
\def\ldots{\DN@{\extrap@\plainldots@}%
 \FN@\next@}
\def\cdots{\DN@{\extrap@\plaincdots@}%
 \FN@\next@}

\def\dotsc{\relaxnext@
 \DN@{\ifx\next;\plainldots@\,\else
  \ifx\next.\plainldots@\,\else\extra@\plainldots@
  \ifextra@\,\fi\fi\fi}%
 \FN@\next@}
\def\cdot{\mathchar"2201 }

\message{special superscripts,}
\def\dddot#1{{\mathop{#1}\limits^{\vbox to-1.4\ex@{\kern-\tw@\ex@
 \hbox{\rm...}\vss}}}}
\def\ddddot#1{{\mathop{#1}\limits^{\vbox to-1.4\ex@{\kern-\tw@\ex@
 \hbox{\rm....}\vss}}}}
\def\sphat{^{\mathchoice{}{}%
 {\,\,\botsmash{\hbox{\lower4\ex@\hbox{$\m@th\widehat{\null}$}}}}%
 {\,\botsmash{\hbox{\lower3\ex@\hbox{$\m@th\hat{\null}$}}}}}}

\def\spacute{^{\!\botsmash{\hbox{\lower\@ne ex\hbox{\'{}}}}}}
\def\spgrave{^{\mathchoice{}{}{}{\!}%
 \botsmash{\hbox{\lower\@ne ex\hbox{\`{}}}}}}
\def\spdot{^{\hbox{\raise\ex@\hbox{\rm.}}}}
\def\spddot{^{\hbox{\raise\ex@\hbox{\rm..}}}}
\def\spdddot{^{\hbox{\raise\ex@\hbox{\rm...}}}}
\def\spddddot{^{\hbox{\raise\ex@\hbox{\rm....}}}}
\def\spbreve{^{\!\botsmash{\hbox{\lower4\ex@\hbox{\u{}}}}}}

\message{\string\text,}
\def\textonlyfont@#1#2{\def#1{\RIfM@
 \Err@{Use \string#1\space only in text}\else#2\fi}}
\textonlyfont@\rm\tenrm
\textonlyfont@\it\tenit
\textonlyfont@\sl\tensl
\textonlyfont@\bf\tenbf
\def\oldnos#1{\RIfM@{\mathcode`\,="013B \fam\@ne#1}\else
 \leavevmode\hbox{$\m@th\mathcode`\,="013B \fam\@ne#1$}\fi}
\def\text{\RIfM@\expandafter\text@\else\expandafter\text@@\fi}
\def\text@@#1{\leavevmode\hbox{#1}}
\def\mathhexbox@#1#2#3{\text{$\m@th\mathchar"#1#2#3$}}
\def\dag{{\mathhexbox@279}}
\def\ddag{{\mathhexbox@27A}}
\def\S{{\mathhexbox@278}}
\def\P{{\mathhexbox@27B}}
\newif\iffirstchoice@
\firstchoice@true
\def\text@#1{\mathchoice
 {\hbox{\everymath{\displaystyle}\def\textfonti{\the\textfont\@ne}%
  \def\textfontii{\the\textfont\tw@}\textdef@@ T#1}}
 {\hbox{\firstchoice@false
  \everymath{\textstyle}\def\textfonti{\the\textfont\@ne}%
  \def\textfontii{\the\textfont\tw@}\textdef@@ T#1}}
 {\hbox{\firstchoice@false
  \everymath{\scriptstyle}\def\textfonti{\the\scriptfont\@ne}%
  \def\textfontii{\the\scriptfont\tw@}\textdef@@ S\rm#1}}
 {\hbox{\firstchoice@false
  \everymath{\scriptscriptstyle}\def\textfonti
  {\the\scriptscriptfont\@ne}%
  \def\textfontii{\the\scriptscriptfont\tw@}\textdef@@ s\rm#1}}}
\def\textdef@@#1{\textdef@#1\rm\textdef@#1\bf\textdef@#1\sl\textdef@#1\it}
\def\rmfam{0}
\def\textdef@#1#2{%
 \DN@{\csname\expandafter\eat@\string#2fam\endcsname}%
 \if S#1\edef#2{\the\scriptfont\next@\relax}%
 \else\if s#1\edef#2{\the\scriptscriptfont\next@\relax}%
 \else\edef#2{\the\textfont\next@\relax}\fi\fi}
\scriptfont\itfam\tenit \scriptscriptfont\itfam\tenit
\scriptfont\slfam\tensl \scriptscriptfont\slfam\tensl
\newif\iftopfolded@
\newif\ifbotfolded@
\def\topfoldedtext{\topfolded@true\botfolded@false\foldedtext@}
\def\botfoldedtext{\botfolded@true\topfolded@false\foldedtext@}
\def\foldedtext{\topfolded@false\botfolded@false\foldedtext@}
\Invalid@\foldedwidth
\def\foldedtext@{\relaxnext@
 \DN@{\ifx\next\foldedwidth\let\next@\nextii@\else
  \DN@{\nextii@\foldedwidth{.3\hsize}}\fi\next@}%
 \DNii@\foldedwidth##1##2{\setbox\z@\vbox
  {\normalbaselines\hsize##1\relax
  \tolerance1600 \noindent\ignorespaces##2}\ifbotfolded@\boxz@\else
  \iftopfolded@\vtop{\unvbox\z@}\else\vcenter{\boxz@}\fi\fi}%
 \FN@\next@}
\message{math font commands,}
\def\bold{\RIfM@\expandafter\bold@\else
 \expandafter\nonmatherr@\expandafter\bold\fi}
\def\bold@#1{{\bold@@{#1}}}
\def\bold@@#1{\fam\bffam\relax#1}
\def\slanted{\RIfM@\expandafter\slanted@\else
 \expandafter\nonmatherr@\expandafter\slanted\fi}
\def\slanted@#1{{\slanted@@{#1}}}
\def\slanted@@#1{\fam\slfam\relax#1}
\def\roman{\RIfM@\expandafter\roman@\else
 \expandafter\nonmatherr@\expandafter\roman\fi}
\def\roman@#1{{\roman@@{#1}}}
\def\roman@@#1{\fam\rmfam\relax#1}
\def\italic{\RIfM@\expandafter\italic@\else
 \expandafter\nonmatherr@\expandafter\italic\fi}
\def\italic@#1{{\italic@@{#1}}}
\def\italic@@#1{\fam\itfam\relax#1}
\def\Cal{\RIfM@\expandafter\Cal@\else
 \expandafter\nonmatherr@\expandafter\Cal\fi}
\def\Cal@#1{{\Cal@@{#1}}}
\def\Cal@@#1{\noaccents@\fam\tw@#1}
\mathchardef\Gamma="0000
\mathchardef\Delta="0001
\mathchardef\Theta="0002
\mathchardef\Lambda="0003
\mathchardef\Xi="0004
\mathchardef\Pi="0005
\mathchardef\Sigma="0006
\mathchardef\Upsilon="0007
\mathchardef\Phi="0008
\mathchardef\Psi="0009
\mathchardef\Omega="000A
\mathchardef\varGamma="0100
\mathchardef\varDelta="0101
\mathchardef\varTheta="0102
\mathchardef\varLambda="0103
\mathchardef\varXi="0104
\mathchardef\varPi="0105
\mathchardef\varSigma="0106
\mathchardef\varUpsilon="0107
\mathchardef\varPhi="0108
\mathchardef\varPsi="0109
\mathchardef\varOmega="010A
\let\alloc@@\alloc@
\def\hexnumber@#1{\ifcase#1 0\or 1\or 2\or 3\or 4\or 5\or 6\or 7\or 8\or
 9\or A\or B\or C\or D\or E\or F\fi}
\def\loadmsam{%
 \font@\tenmsa=msam10
 \font@\sevenmsa=msam7
 \font@\fivemsa=msam5
 \alloc@@8\fam\chardef\sixt@@n\msafam
 \textfont\msafam=\tenmsa
 \scriptfont\msafam=\sevenmsa
 \scriptscriptfont\msafam=\fivemsa
 \edef\next{\hexnumber@\msafam}%
 \mathchardef\dabar@"0\next39
 \edef\dashrightarrow{\mathrel{\dabar@\dabar@\mathchar"0\next4B}}%
 \edef\dashleftarrow{\mathrel{\mathchar"0\next4C\dabar@\dabar@}}%
 \let\dasharrow\dashrightarrow
 \edef\ulcorner{\delimiter"4\next70\next70 }%
 \edef\urcorner{\delimiter"5\next71\next71 }%
 \edef\llcorner{\delimiter"4\next78\next78 }%
 \edef\lrcorner{\delimiter"5\next79\next79 }%
 \edef\yen{{\noexpand\mathhexbox@\next55}}%
 \edef\checkmark{{\noexpand\mathhexbox@\next58}}%
 \edef\circledR{{\noexpand\mathhexbox@\next72}}%
 \edef\maltese{{\noexpand\mathhexbox@\next7A}}%
 \global\let\loadmsam\empty}%
\def\loadmsbm{%
 \font@\tenmsb=msbm10 \font@\sevenmsb=msbm7 \font@\fivemsb=msbm5
 \alloc@@8\fam\chardef\sixt@@n\msbfam
 \textfont\msbfam=\tenmsb
 \scriptfont\msbfam=\sevenmsb \scriptscriptfont\msbfam=\fivemsb
 \global\let\loadmsbm\empty
 }
\def\widehat#1{\ifx\undefined\msbfam \DN@{362}%
  \else \setboxz@h{$\m@th#1$}%
    \edef\next@{\ifdim\wdz@>\tw@ em%
        \hexnumber@\msbfam 5B%
      \else 362\fi}\fi
  \mathaccent"0\next@{#1}}
\def\widetilde#1{\ifx\undefined\msbfam \DN@{365}%
  \else \setboxz@h{$\m@th#1$}%
    \edef\next@{\ifdim\wdz@>\tw@ em%
        \hexnumber@\msbfam 5D%
      \else 365\fi}\fi
  \mathaccent"0\next@{#1}}
\message{\string\newsymbol,}
\def\newsymbol#1#2#3#4#5{\define#1{}%
  \count@#2\relax \advance\count@\m@ne 
 \ifcase\count@
   \ifx\undefined\msafam\loadmsam\fi \let\next@\msafam
 \or \ifx\undefined\msbfam\loadmsbm\fi \let\next@\msbfam
 \else  \Err@{\Invalid@@\string\newsymbol}\let\next@\tw@\fi
 \mathchardef#1="#3\hexnumber@\next@#4#5\space}

\def\Bbb{\RIfM@\expandafter\Bbb@\else
 \expandafter\nonmatherr@\expandafter\Bbb\fi}
\def\Bbb@#1{{\Bbb@@{#1}}}
\def\Bbb@@#1{\noaccents@\fam\msbfam\relax#1}
\message{bold Greek and bold symbols,}
\def\loadbold{%
 \font@\tencmmib=cmmib10 \font@\sevencmmib=cmmib7 \font@\fivecmmib=cmmib5
 \skewchar\tencmmib'177 \skewchar\sevencmmib'177 \skewchar\fivecmmib'177
 \alloc@@8\fam\chardef\sixt@@n\cmmibfam
 \textfont\cmmibfam\tencmmib
 \scriptfont\cmmibfam\sevencmmib \scriptscriptfont\cmmibfam\fivecmmib
 \font@\tencmbsy=cmbsy10 \font@\sevencmbsy=cmbsy7 \font@\fivecmbsy=cmbsy5
 \skewchar\tencmbsy'60 \skewchar\sevencmbsy'60 \skewchar\fivecmbsy'60
 \alloc@@8\fam\chardef\sixt@@n\cmbsyfam
 \textfont\cmbsyfam\tencmbsy
 \scriptfont\cmbsyfam\sevencmbsy \scriptscriptfont\cmbsyfam\fivecmbsy
 \let\loadbold\empty
}
\def\boldnotloaded#1{\Err@{\ifcase#1\or First\else Second\fi
       bold symbol font not loaded}}
\def\mathchari@#1#2#3{\ifx\undefined\cmmibfam
    \boldnotloaded@\@ne
  \else\mathchar"#1\hexnumber@\cmmibfam#2#3\space \fi}
\def\mathcharii@#1#2#3{\ifx\undefined\cmbsyfam
    \boldnotloaded\tw@
  \else \mathchar"#1\hexnumber@\cmbsyfam#2#3\space\fi}
\edef\bffam@{\hexnumber@\bffam}
\def\boldkey#1{\ifcat\noexpand#1A%
  \ifx\undefined\cmmibfam \boldnotloaded\@ne
  \else {\fam\cmmibfam#1}\fi
 \else
 \ifx#1!\mathchar"5\bffam@21 \else
 \ifx#1(\mathchar"4\bffam@28 \else\ifx#1)\mathchar"5\bffam@29 \else
 \ifx#1+\mathchar"2\bffam@2B \else\ifx#1:\mathchar"3\bffam@3A \else
 \ifx#1;\mathchar"6\bffam@3B \else\ifx#1=\mathchar"3\bffam@3D \else
 \ifx#1?\mathchar"5\bffam@3F \else\ifx#1[\mathchar"4\bffam@5B \else
 \ifx#1]\mathchar"5\bffam@5D \else
 \ifx#1,\mathchari@63B \else
 \ifx#1-\mathcharii@200 \else
 \ifx#1.\mathchari@03A \else
 \ifx#1/\mathchari@03D \else
 \ifx#1<\mathchari@33C \else
 \ifx#1>\mathchari@33E \else
 \ifx#1*\mathcharii@203 \else
 \ifx#1|\mathcharii@06A \else
 \ifx#10\bold0\else\ifx#11\bold1\else\ifx#12\bold2\else\ifx#13\bold3\else
 \ifx#14\bold4\else\ifx#15\bold5\else\ifx#16\bold6\else\ifx#17\bold7\else
 \ifx#18\bold8\else\ifx#19\bold9\else
  \Err@{\string\boldkey\space can't be used with #1}%
 \fi\fi\fi\fi\fi\fi\fi\fi\fi\fi\fi\fi\fi\fi\fi
 \fi\fi\fi\fi\fi\fi\fi\fi\fi\fi\fi\fi\fi\fi}
\def\boldsymbol#1{%
 \DN@{\Err@{You can't use \string\boldsymbol\space with \string#1}#1}%
 \ifcat\noexpand#1A%
   \let\next@\relax
   \ifx\undefined\cmmibfam \boldnotloaded\@ne
   \else {\fam\cmmibfam#1}\fi
 \else
  \xdef\meaning@{\meaning#1.........}%
  \expandafter\math@\meaning@\math@
  \ifmath@
   \expandafter\mathch@\meaning@\mathch@
   \ifmathch@
    \expandafter\boldsymbol@@\meaning@\boldsymbol@@
   \fi
  \else
   \expandafter\macro@\meaning@\macro@
   \expandafter\delim@\meaning@\delim@
   \ifdelim@
    \expandafter\delim@@\meaning@\delim@@
   \else
    \boldsymbol@{#1}%
   \fi
  \fi
 \fi
 \next@}
\def\mathhexboxii@#1#2{\ifx\undefined\cmbsyfam
    \boldnotloaded\tw@
  \else \mathhexbox@{\hexnumber@\cmbsyfam}{#1}{#2}\fi}
\def\boldsymbol@#1{\let\next@\relax\let\next#1%
 \ifx\next\cdot\mathcharii@201 \else
 \ifx\next\prime{{\null\mathcharii@030 \null}}\else
 \ifx\next\lbrack\mathchar"4\bffam@5B \else
 \ifx\next\rbrack\mathchar"5\bffam@5D \else
 \ifx\next\{\mathcharii@466 \else
 \ifx\next\lbrace\mathcharii@466 \else
 \ifx\next\}\mathcharii@567 \else
 \ifx\next\rbrace\mathcharii@567 \else
 \ifx\next\surd{{\mathcharii@170}}\else
 \ifx\next\S{{\mathhexboxii@78}}\else
 \ifx\next\P{{\mathhexboxii@7B}}\else
 \ifx\next\dag{{\mathhexboxii@79}}\else
 \ifx\next\ddag{{\mathhexboxii@7A}}\else
 \DN@{\Err@{You can't use \string\boldsymbol\space with \string#1}#1}%
 \fi\fi\fi\fi\fi\fi\fi\fi\fi\fi\fi\fi\fi}
\def\boldsymbol@@#1.#2\boldsymbol@@{\classnum@#1 \count@@@\classnum@        
 \divide\classnum@4096 \count@\classnum@                                    
 \multiply\count@4096 \advance\count@@@-\count@ \count@@\count@@@           
 \divide\count@@@\@cclvi \count@\count@@                                    
 \multiply\count@@@\@cclvi \advance\count@@-\count@@@                       
 \divide\count@@@\@cclvi                                                    
 \multiply\classnum@4096 \advance\classnum@\count@@                         
 \ifnum\count@@@=\z@                                                        
  \count@"\bffam@ \multiply\count@\@cclvi
  \advance\classnum@\count@
  \DN@{\mathchar\number\classnum@}%
 \else
  \ifnum\count@@@=\@ne                                                      
   \ifx\undefined\cmmibfam \DN@{\boldnotloaded\@ne}%
   \else \count@\cmmibfam \multiply\count@\@cclvi
     \advance\classnum@\count@
     \DN@{\mathchar\number\classnum@}\fi
  \else
   \ifnum\count@@@=\tw@                                                    
     \ifx\undefined\cmbsyfam
       \DN@{\boldnotloaded\tw@}%
     \else
       \count@\cmbsyfam \multiply\count@\@cclvi
       \advance\classnum@\count@
       \DN@{\mathchar\number\classnum@}%
     \fi
  \fi
 \fi
\fi}
\newif\ifdelim@
\newcount\delimcount@
{\uccode`6=`\\ \uccode`7=`d \uccode`8=`e \uccode`9=`l
 \uppercase{\gdef\delim@#1#2#3#4#5\delim@
  {\delim@false\ifx 6#1\ifx 7#2\ifx 8#3\ifx 9#4\delim@true
   \xdef\meaning@{#5}\fi\fi\fi\fi}}}
\def\delim@@#1"#2#3#4#5#6\delim@@{\if#32%
\let\next@\relax
 \ifx\undefined\cmbsyfam \boldnotloaded\@ne
 \else \mathcharii@#2#4#5\space \fi\fi}
\def\vert{\delimiter"026A30C }
\def\Vert{\delimiter"026B30D }
\let\|\Vert
\def\backslash{\delimiter"026E30F }
\def\boldkeydots@#1{\bold@true\let\next=#1\let\delayed@=#1\mdots@@
 \boldkey#1\bold@false}  
\def\boldsymboldots@#1{\bold@true\let\next#1\let\delayed@#1\mdots@@
 \boldsymbol#1\bold@false}
\message{Euler fonts,}

\def\frak{\mathfont@\frak}

\def\loadmathfont#1{%
   \expandafter\font@\csname ten#1\endcsname=#110
   \expandafter\font@\csname seven#1\endcsname=#17
   \expandafter\font@\csname five#1\endcsname=#15
   \edef\next{\noexpand\alloc@@8\fam\chardef\sixt@@n
     \expandafter\noexpand\csname#1fam\endcsname}%
   \next
   \textfont\csname#1fam\endcsname \csname ten#1\endcsname
   \scriptfont\csname#1fam\endcsname \csname seven#1\endcsname
   \scriptscriptfont\csname#1fam\endcsname \csname five#1\endcsname
   \expandafter\def\csname #1\expandafter\endcsname\expandafter{%
      \expandafter\mathfont@\csname#1\endcsname}%
 \expandafter\gdef\csname load#1\endcsname{}%
}
\def\mathfont@#1{\RIfM@\expandafter\mathfont@@\expandafter#1\else
  \expandafter\nonmatherr@\expandafter#1\fi}
\def\mathfont@@#1#2{{\mathfont@@@#1{#2}}}
\def\mathfont@@@#1#2{\noaccents@
   \fam\csname\expandafter\eat@\string#1fam\endcsname
   \relax#2}
\message{math accents,}
\def\accentclass@{7}
\def\noaccents@{\def\accentclass@{0}}
\def\makeacc@#1#2{\def#1{\mathaccent"\accentclass@#2 }}
\makeacc@\hat{05E}
\makeacc@\check{014}
\makeacc@\tilde{07E}
\makeacc@\acute{013}
\makeacc@\grave{012}
\makeacc@\dot{05F}
\makeacc@\ddot{07F}
\makeacc@\breve{015}
\makeacc@\bar{016}

\newcount\skewcharcount@
\newcount\familycount@
\def\theskewchar@{\familycount@\@ne
 \global\skewcharcount@\the\skewchar\textfont\@ne                           
 \ifnum\fam>\m@ne\ifnum\fam<16
  \global\familycount@\the\fam\relax
  \global\skewcharcount@\the\skewchar\textfont\the\fam\relax\fi\fi          
 \ifnum\skewcharcount@>\m@ne
  \ifnum\skewcharcount@<128
  \multiply\familycount@256
  \global\advance\skewcharcount@\familycount@
  \global\advance\skewcharcount@28672
  \mathchar\skewcharcount@\else
  \global\skewcharcount@\m@ne\fi\else
 \global\skewcharcount@\m@ne\fi}                                            
\newcount\pointcount@
\def\getpoints@#1.#2\getpoints@{\pointcount@#1 }
\newdimen\accentdimen@
\newcount\accentmu@
\def\dimentomu@{\multiply\accentdimen@ 100
 \expandafter\getpoints@\the\accentdimen@\getpoints@
 \multiply\pointcount@18
 \divide\pointcount@\@m
 \global\accentmu@\pointcount@}
\def\Makeacc@#1#2{\def#1{\RIfM@\DN@{\mathaccent@
 {"\accentclass@#2 }}\else\DN@{\nonmatherr@{#1}}\fi\next@}}
\def\unbracefonts@{\let\Cal@\Cal@@\let\roman@\roman@@\let\bold@\bold@@
 \let\slanted@\slanted@@}
\def\mathaccent@#1#2{\ifnum\fam=\m@ne\xdef\thefam@{1}\else
 \xdef\thefam@{\the\fam}\fi                                                 
 \accentdimen@\z@                                                           
 \setboxz@h{\unbracefonts@$\m@th\fam\thefam@\relax#2$}
 \ifdim\accentdimen@=\z@\DN@{\mathaccent#1{#2}}
  \setbox@ne\hbox{\unbracefonts@$\m@th\fam\thefam@\relax#2\theskewchar@$}
  \setbox\tw@\hbox{$\m@th\ifnum\skewcharcount@=\m@ne\else
   \mathchar\skewcharcount@\fi$}
  \global\accentdimen@\wd@ne\global\advance\accentdimen@-\wdz@
  \global\advance\accentdimen@-\wd\tw@                                     
  \global\multiply\accentdimen@\tw@
  \dimentomu@\global\advance\accentmu@\@ne                                 
 \else\DN@{{\mathaccent#1{#2\mkern\accentmu@ mu}%
    \mkern-\accentmu@ mu}{}}\fi                                             
 \next@}\Makeacc@\Hat{05E}
\Makeacc@\Check{014}
\Makeacc@\Tilde{07E}
\Makeacc@\Acute{013}
\Makeacc@\Grave{012}
\Makeacc@\Dot{05F}
\Makeacc@\Ddot{07F}
\Makeacc@\Breve{015}
\Makeacc@\Bar{016}
\def\Vec{\RIfM@\DN@{\mathaccent@{"017E }}\else
 \DN@{\nonmatherr@\Vec}\fi\next@}
\def\accentedsymbol#1#2{\csname newbox\expandafter\endcsname
  \csname\expandafter\eat@\string#1@box\endcsname
 \expandafter\setbox\csname\expandafter\eat@
  \string#1@box\endcsname\hbox{$\m@th#2$}\define
  #1{\copy\csname\expandafter\eat@\string#1@box\endcsname{}}}
\message{roots,}
\def\sqrt#1{\radical"270370 {#1}}
\let\underline@\underline
\let\overline@\overline
\def\underline#1{\underline@{#1}}
\def\overline#1{\overline@{#1}}
\Invalid@\leftroot
\Invalid@\uproot
\newcount\uproot@
\newcount\leftroot@
\def\root{\relaxnext@
  \DN@{\ifx\next\uproot\let\next@\nextii@\else
   \ifx\next\leftroot\let\next@\nextiii@\else
   \let\next@\plainroot@\fi\fi\next@}%
  \DNii@\uproot##1{\uproot@##1\relax\FN@\nextiv@}%
  \def\nextiv@{\ifx\next\space@\DN@. {\FN@\nextv@}\else
   \DN@.{\FN@\nextv@}\fi\next@.}%
  \def\nextv@{\ifx\next\leftroot\let\next@\nextvi@\else
   \let\next@\plainroot@\fi\next@}%
  \def\nextvi@\leftroot##1{\leftroot@##1\relax\plainroot@}%
   \def\nextiii@\leftroot##1{\leftroot@##1\relax\FN@\nextvii@}%
  \def\nextvii@{\ifx\next\space@
   \DN@. {\FN@\nextviii@}\else
   \DN@.{\FN@\nextviii@}\fi\next@.}%
  \def\nextviii@{\ifx\next\uproot\let\next@\nextix@\else
   \let\next@\plainroot@\fi\next@}%
  \def\nextix@\uproot##1{\uproot@##1\relax\plainroot@}%
  \bgroup\uproot@\z@\leftroot@\z@\FN@\next@}
\def\plainroot@#1\of#2{\setbox\rootbox\hbox{$\m@th\scriptscriptstyle{#1}$}%
 \mathchoice{\r@@t\displaystyle{#2}}{\r@@t\textstyle{#2}}
 {\r@@t\scriptstyle{#2}}{\r@@t\scriptscriptstyle{#2}}\egroup}
\def\r@@t#1#2{\setboxz@h{$\m@th#1\sqrt{#2}$}%
 \dimen@\ht\z@\advance\dimen@-\dp\z@
 \setbox@ne\hbox{$\m@th#1\mskip\uproot@ mu$}\advance\dimen@ 1.667\wd@ne
 \mkern-\leftroot@ mu\mkern5mu\raise.6\dimen@\copy\rootbox
 \mkern-10mu\mkern\leftroot@ mu\boxz@}
\def\boxed#1{\setboxz@h{$\m@th\displaystyle{#1}$}\dimen@.4\ex@
 \advance\dimen@3\ex@\advance\dimen@\dp\z@
 \hbox{\lower\dimen@\hbox{%
 \vbox{\hrule height.4\ex@
 \hbox{\vrule width.4\ex@\hskip3\ex@\vbox{\vskip3\ex@\boxz@\vskip3\ex@}%
 \hskip3\ex@\vrule width.4\ex@}\hrule height.4\ex@}%
 }}}
\message{commutative diagrams,}
\let\ampersand@\relax
\newdimen\minaw@
\minaw@11.11128\ex@
\newdimen\minCDaw@
\minCDaw@2.5pc
\def\minCDarrowwidth#1{\RIfMIfI@\onlydmatherr@\minCDarrowwidth
 \else\minCDaw@#1\relax\fi\else\onlydmatherr@\minCDarrowwidth\fi}
\newif\ifCD@
\def\CD{\bgroup\vspace@\relax\let\ampersand@&\iffalse}\fi
 \CD@true\vcenter\bgroup\Let@\tabskip\z@skip\baselineskip20\ex@
 \lineskip3\ex@\lineskiplimit3\ex@\halign\bgroup
 &\hfill$\m@th##$\hfill\crcr}
\def\endCD{\crcr\egroup\egroup\egroup}
\newdimen\bigaw@
\atdef@>#1>#2>{\ampersand@                                                  
 \setboxz@h{$\m@th\ssize\;{#1}\;\;$}
 \setbox@ne\hbox{$\m@th\ssize\;{#2}\;\;$}
 \setbox\tw@\hbox{$\m@th#2$}
 \ifCD@\global\bigaw@\minCDaw@\else\global\bigaw@\minaw@\fi                 
 \ifdim\wdz@>\bigaw@\global\bigaw@\wdz@\fi
 \ifdim\wd@ne>\bigaw@\global\bigaw@\wd@ne\fi                                
 \ifCD@\enskip\fi                                                           
 \ifdim\wd\tw@>\z@
  \mathrel{\mathop{\hbox to\bigaw@{\rightarrowfill@\displaystyle}}%
    \limits^{#1}_{#2}}
 \else\mathrel{\mathop{\hbox to\bigaw@{\rightarrowfill@\displaystyle}}%
    \limits^{#1}}\fi                                                        
 \ifCD@\enskip\fi                                                          
 \ampersand@}                                                              
\atdef@<#1<#2<{\ampersand@\setboxz@h{$\m@th\ssize\;\;{#1}\;$}%
 \setbox@ne\hbox{$\m@th\ssize\;\;{#2}\;$}\setbox\tw@\hbox{$\m@th#2$}%
 \ifCD@\global\bigaw@\minCDaw@\else\global\bigaw@\minaw@\fi
 \ifdim\wdz@>\bigaw@\global\bigaw@\wdz@\fi
 \ifdim\wd@ne>\bigaw@\global\bigaw@\wd@ne\fi
 \ifCD@\enskip\fi
 \ifdim\wd\tw@>\z@
  \mathrel{\mathop{\hbox to\bigaw@{\leftarrowfill@\displaystyle}}%
       \limits^{#1}_{#2}}\else
  \mathrel{\mathop{\hbox to\bigaw@{\leftarrowfill@\displaystyle}}%
       \limits^{#1}}\fi
 \ifCD@\enskip\fi\ampersand@}
\begingroup
 \catcode`\~=\active \lccode`\~=`\@
 \lowercase{%
  \global\atdef@)#1)#2){~>#1>#2>}
  \global\atdef@(#1(#2({~<#1<#2<}}
\endgroup
\atdef@ A#1A#2A{\llap{$\m@th\vcenter{\hbox
 {$\ssize#1$}}$}\Big\uparrow\rlap{$\m@th\vcenter{\hbox{$\ssize#2$}}$}&&}
\atdef@ V#1V#2V{\llap{$\m@th\vcenter{\hbox
 {$\ssize#1$}}$}\Big\downarrow\rlap{$\m@th\vcenter{\hbox{$\ssize#2$}}$}&&}
\atdef@={&\enskip\mathrel
 {\vbox{\hrule width\minCDaw@\vskip3\ex@\hrule width
 \minCDaw@}}\enskip&}
\atdef@|{\Big\Vert&&}
\atdef@\vert{\Big\Vert&&}
\def\pretend#1\haswidth#2{\setboxz@h{$\m@th\scriptstyle{#2}$}\hbox
 to\wdz@{\hfill$\m@th\scriptstyle{#1}$\hfill}}
\message{poor man's bold,}
\def\pmb{\RIfM@\expandafter\mathpalette\expandafter\pmb@\else
 \expandafter\pmb@@\fi}
\def\pmb@@#1{\leavevmode\setboxz@h{#1}%
   \dimen@-\wdz@
   \kern-.5\ex@\copy\z@
   \kern\dimen@\kern.25\ex@\raise.4\ex@\copy\z@
   \kern\dimen@\kern.25\ex@\box\z@
}
\def\binrel@@#1{\ifdim\wd2<\z@\mathbin{#1}\else\ifdim\wd\tw@>\z@
 \mathrel{#1}\else{#1}\fi\fi}
\newdimen\pmbraise@
\def\pmb@#1#2{\setbox\thr@@\hbox{$\m@th#1{#2}$}%
 \setbox4\hbox{$\m@th#1\mkern.5mu$}\pmbraise@\wd4\relax
 \binrel@{#2}%
 \dimen@-\wd\thr@@
   \binrel@@{%
   \mkern-.8mu\copy\thr@@
   \kern\dimen@\mkern.4mu\raise\pmbraise@\copy\thr@@
   \kern\dimen@\mkern.4mu\box\thr@@
}}
\def\documentstyle#1{\W@{}\input #1.sty\relax}
\message{syntax check,}
\font\dummyft@=dummy
\fontdimen1 \dummyft@=\z@
\fontdimen2 \dummyft@=\z@
\fontdimen3 \dummyft@=\z@
\fontdimen4 \dummyft@=\z@
\fontdimen5 \dummyft@=\z@
\fontdimen6 \dummyft@=\z@
\fontdimen7 \dummyft@=\z@
\fontdimen8 \dummyft@=\z@
\fontdimen9 \dummyft@=\z@
\fontdimen10 \dummyft@=\z@
\fontdimen11 \dummyft@=\z@
\fontdimen12 \dummyft@=\z@
\fontdimen13 \dummyft@=\z@
\fontdimen14 \dummyft@=\z@
\fontdimen15 \dummyft@=\z@
\fontdimen16 \dummyft@=\z@
\fontdimen17 \dummyft@=\z@
\fontdimen18 \dummyft@=\z@
\fontdimen19 \dummyft@=\z@
\fontdimen20 \dummyft@=\z@
\fontdimen21 \dummyft@=\z@
\fontdimen22 \dummyft@=\z@
\def\fontlist@{\\{\tenrm}\\{\sevenrm}\\{\fiverm}\\{\teni}\\{\seveni}%
 \\{\fivei}\\{\tensy}\\{\sevensy}\\{\fivesy}\\{\tenex}\\{\tenbf}\\{\sevenbf}%
 \\{\fivebf}\\{\tensl}\\{\tenit}}
\def\font@#1=#2 {\rightappend@#1\to\fontlist@\font#1=#2 }
\def\dodummy@{{\def\\##1{\global\let##1\dummyft@}\fontlist@}}
\def\nopages@{\output{\setbox\z@\box\@cclv \deadcycles\z@}%
 \alloc@5\toks\toksdef\@cclvi\output}
\let\galleys\nopages@
\newif\ifsyntax@
\newcount\countxviii@
\def\syntax{\syntax@true\dodummy@\countxviii@\count18
 \loop\ifnum\countxviii@>\m@ne\textfont\countxviii@=\dummyft@
 \scriptfont\countxviii@=\dummyft@\scriptscriptfont\countxviii@=\dummyft@
 \advance\countxviii@\m@ne\repeat                                           
 \dummyft@\tracinglostchars\z@\nopages@\frenchspacing\hbadness\@M}
\def\first@#1#2\end{#1}
\def\printoptions{\W@{Do you want S(yntax check),
  G(alleys) or P(ages)?}%
 \message{Type S, G or P, followed by <return>: }%
 \begingroup 
 \endlinechar\m@ne 
 \read\m@ne to\ans@
 \edef\ans@{\uppercase{\def\noexpand\ans@{%
   \expandafter\first@\ans@ P\end}}}%
 \expandafter\endgroup\ans@
 \if\ans@ P
 \else \if\ans@ S\syntax
 \else \if\ans@ G\galleys
 \else\message{? Unknown option: \ans@; using the `pages' option.}%
 \fi\fi\fi}
\def\alloc@#1#2#3#4#5{\global\advance\count1#1by\@ne
 \ch@ck#1#4#2\allocationnumber=\count1#1
 \global#3#5=\allocationnumber
 \ifalloc@\wlog{\string#5=\string#2\the\allocationnumber}\fi}
\def\document{\def\alloclist@{}\def\fontlist@{}}
\let\enddocument\bye

\let\proclaim\undefined
\let\footnote\undefined
\let\=\undefined
\let\>\undefined

\catcode`\@=\active
\message{... finished}

\expandafter\ifx\csname mathdefs.tex\endcsname\relax
  \expandafter\gdef\csname mathdefs.tex\endcsname{}
\else \message{Hey!  Apparently you were trying to
  \string\input{mathdefs.tex} twice.   This does not make sense.} 
\errmessage{Please edit your file (probably \jobname.tex) and remove
any duplicate ``\string\input'' lines}\endinput\fi




\catcode`\X=12\catcode`\@=11

\def\n@wcount{\alloc@0\count\countdef\insc@unt}
\def\n@wwrite{\alloc@7\write\chardef\sixt@@n}
\def\n@wread{\alloc@6\read\chardef\sixt@@n}
\def\r@s@t{\relax}\def\v@idline{\par}\def\@mputate#1/{#1}
\def\l@c@l#1X{\firstpart.#1}\def\gl@b@l#1X{#1}\def\t@d@l#1X{{}}

\def\crossrefs#1{\ifx\all#1\let\tr@ce=\all\else\def\tr@ce{#1,}\fi
   \n@wwrite\cit@tionsout\openout\cit@tionsout=\jobname.cit 
   \write\cit@tionsout{\tr@ce}\expandafter\setfl@gs\tr@ce,}
\def\setfl@gs#1,{\def\@{#1}\ifx\@\empty\let\next=\relax
   \else\let\next=\setfl@gs\expandafter\xdef
   \csname#1tr@cetrue\endcsname{}\fi\next}
\def\m@ketag#1#2{\expandafter\n@wcount\csname#2tagno\endcsname
     \csname#2tagno\endcsname=0\let\tail=\all\xdef\all{\tail#2,}
   \ifx#1\l@c@l\let\tail=\r@s@t\xdef\r@s@t{\csname#2tagno\endcsname=0\tail}\fi
   \expandafter\gdef\csname#2cite\endcsname##1{\expandafter
     \ifx\csname#2tag##1\endcsname\relax?\else\csname#2tag##1\endcsname\fi
     \expandafter\ifx\csname#2tr@cetrue\endcsname\relax\else
     \write\cit@tionsout{#2tag ##1 cited on page \folio.}\fi}
   \expandafter\gdef\csname#2page\endcsname##1{\expandafter
     \ifx\csname#2page##1\endcsname\relax?\else\csname#2page##1\endcsname\fi
     \expandafter\ifx\csname#2tr@cetrue\endcsname\relax\else
     \write\cit@tionsout{#2tag ##1 cited on page \folio.}\fi}
   \expandafter\gdef\csname#2tag\endcsname##1{\expandafter
      \ifx\csname#2check##1\endcsname\relax
      \expandafter\xdef\csname#2check##1\endcsname{}%
      \else\immediate\write16{Warning: #2tag ##1 used more than once.}\fi
      \multit@g{#1}{#2}##1/X%
      \write\t@gsout{#2tag ##1 assigned number \csname#2tag##1\endcsname\space
      on page \number\count0.}%
   \csname#2tag##1\endcsname}}

\def\multit@g#1#2#3/#4X{\def\t@mp{#4}\ifx\t@mp\empty%
      \global\advance\csname#2tagno\endcsname by 1 
      \expandafter\xdef\csname#2tag#3\endcsname
      {#1\number\csname#2tagno\endcsnameX}%
   \else\expandafter\ifx\csname#2last#3\endcsname\relax
      \expandafter\n@wcount\csname#2last#3\endcsname
      \global\advance\csname#2tagno\endcsname by 1 
      \expandafter\xdef\csname#2tag#3\endcsname
      {#1\number\csname#2tagno\endcsnameX}
      \write\t@gsout{#2tag #3 assigned number \csname#2tag#3\endcsname\space
      on page \number\count0.}\fi
   \global\advance\csname#2last#3\endcsname by 1
   \def\t@mp{\expandafter\xdef\csname#2tag#3/}%
   \expandafter\t@mp\@mputate#4\endcsname
   {\csname#2tag#3\endcsname\lastpart{\csname#2last#3\endcsname}}\fi}
\def\t@gs#1{\def\all{}\m@ketag#1e\m@ketag#1s\m@ketag\t@d@l p
\let\realscite\scite
\let\realstag\stag
   \m@ketag\gl@b@l r \n@wread\t@gsin
   \openin\t@gsin=\jobname.tgs \re@der \closein\t@gsin
   \n@wwrite\t@gsout\openout\t@gsout=\jobname.tgs }
\outer\def\localtags{\t@gs\l@c@l}
\outer\def\globaltags{\t@gs\gl@b@l}
\outer\def\newlocaltag#1{\m@ketag\l@c@l{#1}}
\outer\def\newglobaltag#1{\m@ketag\gl@b@l{#1}}

\newif\ifpr@ 
\def\m@kecs #1tag #2 assigned number #3 on page #4.%
   {\expandafter\gdef\csname#1tag#2\endcsname{#3}
   \expandafter\gdef\csname#1page#2\endcsname{#4}
   \ifpr@\expandafter\xdef\csname#1check#2\endcsname{}\fi}
\def\re@der{\ifeof\t@gsin\let\next=\relax\else
   \read\t@gsin to\t@gline\ifx\t@gline\v@idline\else
   \expandafter\m@kecs \t@gline\fi\let \next=\re@der\fi\next}
\def\pretags#1{\pr@true\pret@gs#1,,}
\def\pret@gs#1,{\def\@{#1}\ifx\@\empty\let\n@xtfile=\relax
   \else\let\n@xtfile=\pret@gs \openin\t@gsin=#1.tgs \message{#1} \re@der 
   \closein\t@gsin\fi \n@xtfile}

\newcount\sectno\sectno=0\newcount\subsectno\subsectno=0
\newif\ifultr@local \def\ultralocal{\ultr@localtrue}
\def\firstpart{\number\sectno}
\def\lastpart#1{\ifcase#1 \or a\or b\or c\or d\or e\or f\or g\or h\or 
   i\or k\or l\or m\or n\or o\or p\or q\or r\or s\or t\or u\or v\or w\or 
   x\or y\or z \fi}

\def\resetall{\global\advance\sectno by 1\subsectno=0
   \gdef\firstpart{\number\sectno}\r@s@t}
\def\resetsub{\global\advance\subsectno by 1
   \gdef\firstpart{\number\sectno.\number\subsectno}\r@s@t}
\def\newsection#1\par{\resetall\vskip0pt plus.3\vsize\penalty-250
   \vskip0pt plus-.3\vsize\bigskip\bigskip
   \message{#1}\leftline{\bf#1}\nobreak\bigskip}
\def\subsection#1\par{\ifultr@local\resetsub\fi
   \vskip0pt plus.2\vsize\penalty-250\vskip0pt plus-.2\vsize
   \bigskip\smallskip\message{#1}\leftline{\bf#1}\nobreak\medskip}


\newdimen\marginshift

\newdimen\margindelta
\newdimen\marginmax
\newdimen\marginmin

\def\margininit{       
\marginmax=3 true cm                  
				      
\margindelta=0.1 true cm              
\marginmin=0.1true cm                 
\marginshift=\marginmin
}    

\def\t@gsjj#1,{\def\@{#1}\ifx\@\empty\let\next=\relax\else\let\next=\t@gsjj
   \def\@@{p}\ifx\@\@@\else
   \expandafter\gdef\csname#1cite\endcsname##1{\citejj{##1}}
   \expandafter\gdef\csname#1page\endcsname##1{?}
   \expandafter\gdef\csname#1tag\endcsname##1{\tagjj{##1}}\fi\fi\next}
\newif\ifshowstuffinmargin
\showstuffinmarginfalse
\def\jjtags{\ifx\shlhetal\relax 
  \else
\ifx\shlhetal\undefinedcontrolseq
\else
\showstuffinmargintrue
\ifx\all\relax\else\expandafter\t@gsjj\all,\fi\fi \fi
}

\def\tagjj#1{\realstag{#1}\mginpar{\zeigen{#1}}}
\def\citejj#1{\rechnen{#1}\mginpar{\zeigen{#1}}}     

\def\rechnen#1{\expandafter\ifx\csname stag#1\endcsname\relax ??\else
                           \csname stag#1\endcsname\fi}

\newdimen\theight

\def\marginfont{\sevenrm}

\def\trymarginbox#1{\setbox0=\hbox{\marginfont\hskip\marginshift #1}%
		\global\marginshift\wd0 
		\global\advance\marginshift\margindelta}

\def \mginpar#1{%
\ifvmode\setbox0\hbox to \hsize{\hfill\rlap{\marginfont\quad#1}}%
\ht0 0cm
\dp0 0cm
\box0\vskip-\baselineskip
\else 
             \vadjust{\trymarginbox{#1}%
		\ifdim\marginshift>\marginmax \global\marginshift\marginmin
			\trymarginbox{#1}%
                \fi
             \theight=\ht0
             \advance\theight by \dp0    \advance\theight by \lineskip
             \kern -\theight \vbox to \theight{\rightline{\rlap{\box0}}%
\vss}}\fi}


\def\t@gsoff#1,{\def\@{#1}\ifx\@\empty\let\next=\relax\else\let\next=\t@gsoff
   \def\@@{p}\ifx\@\@@\else
   \expandafter\gdef\csname#1cite\endcsname##1{\zeigen{##1}}
   \expandafter\gdef\csname#1page\endcsname##1{?}
   \expandafter\gdef\csname#1tag\endcsname##1{\zeigen{##1}}\fi\fi\next}
\def\verbatimtags{\showstuffinmarginfalse
\ifx\all\relax\else\expandafter\t@gsoff\all,\fi}
\def\zeigen#1{\hbox{$\langle$}#1\hbox{$\rangle$}}

\def\margincite#1{\ifshowstuffinmargin\mginpar{\zeigen{#1}}\fi}

\def\margintag#1{\ifshowstuffinmargin\mginpar{\zeigen{#1}}\fi}

\def\(#1){\edef\dot@g{\ifmmode\ifinner(\hbox{\noexpand\etag{#1}})
   \else\noexpand\eqno(\hbox{\noexpand\etag{#1}})\fi
   \else(\noexpand\ecite{#1})\fi}\dot@g}

\newif\ifbr@ck
\def\eat#1{}
\def\[#1]{\br@cktrue[\br@cket#1'X]}
\def\br@cket#1'#2X{\def\temp{#2}\ifx\temp\empty\let\next\eat
   \else\let\next\br@cket\fi
   \ifbr@ck\br@ckfalse\br@ck@t#1,X\else\br@cktrue#1\fi\next#2X}
\def\br@ck@t#1,#2X{\def\temp{#2}\ifx\temp\empty\let\neext\eat
   \else\let\neext\br@ck@t\def\temp{,}\fi
   \def\teemp{#1}\ifx\teemp\empty\else\rcite{#1}\fi\temp\neext#2X}
\def\resetbr@cket{\gdef\[##1]{[\rtag{##1}]}}
\def\references{\resetbr@cket\newsection References\par}

\newtoks\symb@ls\newtoks\s@mb@ls\newtoks\p@gelist\n@wcount\ftn@mber
    \ftn@mber=1\newif\ifftn@mbers\ftn@mbersfalse\newif\ifbyp@ge\byp@gefalse
\def\defm@rk{\ifftn@mbers\n@mberm@rk\else\symb@lm@rk\fi}
\def\n@mberm@rk{\xdef\m@rk{{\the\ftn@mber}}%
    \global\advance\ftn@mber by 1 }
\def\rot@te#1{\let\temp=#1\global#1=\expandafter\r@t@te\the\temp,X}
\def\r@t@te#1,#2X{{#2#1}\xdef\m@rk{{#1}}}
\def\b@@st#1{{$^{#1}$}}\def\str@p#1{#1}
\def\symb@lm@rk{\ifbyp@ge\rot@te\p@gelist\ifnum\expandafter\str@p\m@rk=1 
    \s@mb@ls=\symb@ls\fi\write\f@nsout{\number\count0}\fi \rot@te\s@mb@ls}
\def\byp@ge{\byp@getrue\n@wwrite\f@nsin\openin\f@nsin=\jobname.fns 
    \n@wcount\currentp@ge\currentp@ge=0\p@gelist={0}
    \re@dfns\closein\f@nsin\rot@te\p@gelist
    \n@wread\f@nsout\openout\f@nsout=\jobname.fns }
\def\m@kelist#1X#2{{#1,#2}}
\def\re@dfns{\ifeof\f@nsin\let\next=\relax\else\read\f@nsin to \f@nline
    \ifx\f@nline\v@idline\else\let\t@mplist=\p@gelist
    \ifnum\currentp@ge=\f@nline
    \global\p@gelist=\expandafter\m@kelist\the\t@mplistX0
    \else\currentp@ge=\f@nline
    \global\p@gelist=\expandafter\m@kelist\the\t@mplistX1\fi\fi
    \let\next=\re@dfns\fi\next}
\def\symbols#1{\symb@ls={#1}\s@mb@ls=\symb@ls} 
\def\bigsymbol{\textstyle}
\symbols{\bigsymbol\ast,\dagger,\ddagger,\sharp,\flat,\natural,\star}
\def\ftnumbers{\ftn@mberstrue} \def\ftsymbols{\ftn@mbersfalse}
\def\paginal{\byp@ge} \def\resetftnumbers{\ftn@mber=1}
\def\ftnote#1{\defm@rk\expandafter\expandafter\expandafter\footnote
    \expandafter\b@@st\m@rk{#1}}

\long\def\jump#1\endjump{}
\def\ssum{\mathop{\lower .1em\hbox{$\textstyle\Sigma$}}\nolimits}

\def\qed{\nobreak\kern 1em \vrule height .5em width .5em depth 0em}
\def\newneq{\hbox{\rlap{\hbox to 1\wd9{\hss$=$\hss}}\raise .1em 
   \hbox to 1\wd9{\hss$\scriptscriptstyle/$\hss}}}
\def\subsetne{\setbox9 = \hbox{$\subset$}\mathrel{\hbox{\rlap
   {\lower .4em \newneq}\raise .13em \hbox{$\subset$}}}}
\def\supsetne{\setbox9 = \hbox{$\subset$}\mathrel{\hbox{\rlap
   {\lower .4em \newneq}\raise .13em \hbox{$\supset$}}}}

\def\vbar{\mathchoice{\vrule height6.3ptdepth-.5ptwidth.8pt\kern-.8pt}
   {\vrule height6.3ptdepth-.5ptwidth.8pt\kern-.8pt}
   {\vrule height4.1ptdepth-.35ptwidth.6pt\kern-.6pt}
   {\vrule height3.1ptdepth-.25ptwidth.5pt\kern-.5pt}}
\def\f@dge{\mathchoice{}{}{\mkern.5mu}{\mkern.8mu}}
\def\b@c#1#2{{\rm \mkern#2mu\vbar\mkern-#2mu#1}}
\def\b@b#1{{\rm I\mkern-3.5mu #1}}
\def\b@a#1#2{{\rm #1\mkern-#2mu\f@dge #1}}
\def\bb#1{{\count4=`#1 \advance\count4by-64 \ifcase\count4\or\b@a A{11.5}\or
   \b@b B\or\b@c C{5}\or\b@b D\or\b@b E\or\b@b F \or\b@c G{5}\or\b@b H\or
   \b@b I\or\b@c J{3}\or\b@b K\or\b@b L \or\b@b M\or\b@b N\or\b@c O{5} \or
   \b@b P\or\b@c Q{5}\or\b@b R\or\b@a S{8}\or\b@a T{10.5}\or\b@c U{5}\or
   \b@a V{12}\or\b@a W{16.5}\or\b@a X{11}\or\b@a Y{11.7}\or\b@a Z{7.5}\fi}}

\catcode`\X=11 \catcode`\@=12




\let\thischap\jobname

\def\partof#1{\csname returnthe#1part\endcsname}
\def\chapof#1{\csname returnthe#1chap\endcsname}

\def\setchapter#1,#2,#3.{%
  \expandafter\def\csname returnthe#1part\endcsname{#2}%
  \expandafter\def\csname returnthe#1chap\endcsname{#3}%
}

\setchapter 300a,A,I.
\setchapter 300b,A,II.
\setchapter 300c,A,III.
\setchapter 300d,A,IV.
\setchapter 300e,A,V.
\setchapter 300f,A,VI.
\setchapter 300g,A,VII.
\setchapter   88,B,I.
\setchapter  600,B,II.
\setchapter  705,B,III.

\def\cprefix#1{
\edef\theotherpart{\partof{#1}}\edef\theotherchap{\chapof{#1}}%
\ifx\theotherpart\thispart
   \ifx\theotherchap\thischap 
    \else 
     \theotherchap%
    \fi
   \else 
     \theotherpart.\theotherchap\fi}

\def\sectioncite[#1]#2{%
     \cprefix{#2}#1}

\edef\thispart{\partof{\thischap}}
\edef\thischap{\chapof{\thischap}}


\def\spuriousreset{}


\expandafter\ifx\csname citeadd.tex\endcsname\relax
\expandafter\gdef\csname citeadd.tex\endcsname{}
\else \message{Hey!  Apparently you were trying to
\string\input{citeadd.tex} twice.   This does not make sense.} 
\errmessage{Please edit your file (probably \jobname.tex) and remove
any duplicate ``\string\input'' lines}\endinput\fi

\sectno=-1   
\localtags
\jjtags
\NoBlackBoxes
\define\mr{\medskip\roster}
\define\sn{\smallskip\noindent}
\define\mn{\medskip\noindent}
\define\bn{\bigskip\noindent}
\define\ub{\underbar}
\define\wilog{\text{without loss of generality}}
\define\ermn{\endroster\medskip\noindent}

\define\dbcu{\dsize\bigcup}
\define \nl{\newline}
\magnification=\magstep 1
\documentstyle{amsppt}

{    
\catcode`@11

\ifx\alicetwothousandloaded@\relax
  \endinput\else\global\let\alicetwothousandloaded@\relax\fi

\gdef\subjclass{\let\savedef@\subjclass
 \def\subjclass##1\endsubjclass{\let\subjclass\savedef@
   \toks@{\def\usualspace{{\rm\enspace}}\eightpoint}%
   \toks@@{##1\unskip.}%
   \edef\thesubjclass@{\the\toks@
     \frills@{{\noexpand\rm2000 {\noexpand\it Mathematics Subject
       Classification}.\noexpand\enspace}}%
     \the\toks@@}}%
  \nofrillscheck\subjclass}
} 


\expandafter\ifx\csname alice2jlem.tex\endcsname\relax
  \expandafter\xdef\csname alice2jlem.tex\endcsname{\the\catcode`@}
\else \message{Hey!  Apparently you were trying to
\string\input{alice2jlem.tex}  twice.   This does not make sense.}
\errmessage{Please edit your file (probably \jobname.tex) and remove
any duplicate ``\string\input'' lines}\endinput\fi

\expandafter\ifx\csname bib4plain.tex\endcsname\relax
  \expandafter\gdef\csname bib4plain.tex\endcsname{}
\else \message{Hey!  Apparently you were trying to \string\input
  bib4plain.tex twice.   This does not make sense.}
\errmessage{Please edit your file (probably \jobname.tex) and remove
any duplicate ``\string\input'' lines}\endinput\fi

\def\renewcommand{\newcommand}	       
\edef\cite{\the\catcode`@}%
\catcode`@ = 11
\let\@oldatcatcode = \cite
\chardef\@letter = 11
\chardef\@other = 12
%
%
%
%
\def\@innerdef#1#2{\edef#1{\expandafter\noexpand\csname #2\endcsname}}%
%
%
\@innerdef\@innernewcount{newcount}%
\@innerdef\@innernewdimen{newdimen}%
\@innerdef\@innernewif{newif}%
\@innerdef\@innernewwrite{newwrite}%
%
%
%
\def\@gobble#1{}%
%
%
%
\ifx\inputlineno\@undefined
   \let\@linenumber = \empty 
\else
   \def\@linenumber{\the\inputlineno:\space}%
\fi
%
%
%
\def\@futurenonspacelet#1{\def\cs{#1}%
   \afterassignment\@stepone\let\@nexttoken=
}%
\begingroup 
\def\\{\global\let\@stoken= }%
\\ 
\endgroup
\def\@stepone{\expandafter\futurelet\cs\@steptwo}%
\def\@steptwo{\expandafter\ifx\cs\@stoken\let\@@next=\@stepthree
   \else\let\@@next=\@nexttoken\fi \@@next}%
\def\@stepthree{\afterassignment\@stepone\let\@@next= }%
%
%
%
\def\@getoptionalarg#1{%
   \let\@optionaltemp = #1%
   \let\@optionalnext = \relax
   \@futurenonspacelet\@optionalnext\@bracketcheck
}%
%
%
\def\@bracketcheck{%
   \ifx [\@optionalnext
      \expandafter\@@getoptionalarg
   \else
      \let\@optionalarg = \empty
      \expandafter\@optionaltemp
   \fi
}%
\def\@@getoptionalarg[#1]{%
   \def\@optionalarg{#1}%
   \@optionaltemp
}%
%
%
%
\def\@nnil{\@nil}%
\def\@fornoop#1\@@#2#3{}%
\def\@for#1:=#2\do#3{%
   \edef\@fortmp{#2}%
   \ifx\@fortmp\empty \else
      \expandafter\@forloop#2,\@nil,\@nil\@@#1{#3}%
   \fi
}%
\def\@forloop#1,#2,#3\@@#4#5{\def#4{#1}\ifx #4\@nnil \else
       #5\def#4{#2}\ifx #4\@nnil \else#5\@iforloop #3\@@#4{#5}\fi\fi
}%
\def\@iforloop#1,#2\@@#3#4{\def#3{#1}\ifx #3\@nnil
       \let\@nextwhile=\@fornoop \else
      #4\relax\let\@nextwhile=\@iforloop\fi\@nextwhile#2\@@#3{#4}%
}%
%
%
%
\@innernewif\if@fileexists
\def\@testfileexistence{\@getoptionalarg\@finishtestfileexistence}%
\def\@finishtestfileexistence#1{%
   \begingroup
      \def\extension{#1}%
      \immediate\openin0 =
         \ifx\@optionalarg\empty\jobname\else\@optionalarg\fi
         \ifx\extension\empty \else .#1\fi
         \space
      \ifeof 0
         \global\@fileexistsfalse
      \else
         \global\@fileexiststrue
      \fi
      \immediate\closein0
   \endgroup
}%
%
%
%
%
\def\bibliographystyle#1{%
   \@readauxfile
   \@writeaux{\string\bibstyle{#1}}%
}%
\let\bibstyle = \@gobble
%
%
\let\bblfilebasename = \jobname
\def\bibliography#1{%
   \@readauxfile
   \@writeaux{\string\bibdata{#1}}%
   \@testfileexistence[\bblfilebasename]{bbl}%
   \if@fileexists
      \nobreak
      \@readbblfile
   \fi
}%
\let\bibdata = \@gobble
%
%
\def\nocite#1{%
   \@readauxfile
   \@writeaux{\string\citation{#1}}%
}%
\@innernewif\if@notfirstcitation
%
%
\def\cite{\@getoptionalarg\@cite}%
%
%
\def\@cite#1{%
   \let\@citenotetext = \@optionalarg
   \printcitestart
   \nocite{#1}%
   \@notfirstcitationfalse
   \@for \@citation :=#1\do
   {%
      \expandafter\@onecitation\@citation\@@
   }%
   \ifx\empty\@citenotetext\else
      \printcitenote{\@citenotetext}%
   \fi
   \printcitefinish
}%
\newif\ifweareinprivate
\weareinprivatetrue
\ifx\shlhetal\undefinedcontrolseq\weareinprivatefalse\fi
\ifx\shlhetal\relax\weareinprivatefalse\fi
\def\@onecitation#1\@@{%
   \if@notfirstcitation
      \printbetweencitations
   \fi
   \expandafter \ifx \csname\@citelabel{#1}\endcsname \relax
      \if@citewarning
         \message{\@linenumber Undefined citation `#1'.}%
      \fi
     \ifweareinprivate
      \expandafter\gdef\csname\@citelabel{#1}\endcsname{%
\strut 
\vadjust{\vskip-\dp\strutbox
\vbox to 0pt{\vss\parindent0cm \leftskip=\hsize 
\advance\leftskip3mm
\advance\hsize 4cm\strut\openup-4pt 
\rightskip 0cm plus 1cm minus 0.5cm ?  #1 ?\strut}}
         {\tt
            \escapechar = -1
            \nobreak\hskip0pt\pfeilsw
            \expandafter\string\csname#1\endcsname
             \pfeilso
            \nobreak\hskip0pt
         }%
      }%
     \else  
      \expandafter\gdef\csname\@citelabel{#1}\endcsname{%
            {\tt\expandafter\string\csname#1\endcsname}
      }%
     \fi  
   \fi
   \csname\@citelabel{#1}\endcsname
   \@notfirstcitationtrue
}%
%
%
\def\@citelabel#1{b@#1}%
%
%
\def\@citedef#1#2{\expandafter\gdef\csname\@citelabel{#1}\endcsname{#2}}%
%
%
%
\def\@readbblfile{%
   \ifx\@itemnum\@undefined
      \@innernewcount\@itemnum
   \fi
   \begingroup
      \def\begin##1##2{%
         \setbox0 = \hbox{\biblabelcontents{##2}}%
         \biblabelwidth = \wd0
      }%
      \def\end##1{}
      %
      %
      \@itemnum = 0
      \def\bibitem{\@getoptionalarg\@bibitem}%
      \def\@bibitem{%
         \ifx\@optionalarg\empty
            \expandafter\@numberedbibitem
         \else
            \expandafter\@alphabibitem
         \fi
      }%
      \def\@alphabibitem##1{%
         \expandafter \xdef\csname\@citelabel{##1}\endcsname {\@optionalarg}%
         \ifx\biblabelprecontents\@undefined
            \let\biblabelprecontents = \relax
         \fi
         \ifx\biblabelpostcontents\@undefined
            \let\biblabelpostcontents = \hss
         \fi
         \@finishbibitem{##1}%
      }%
      \def\@numberedbibitem##1{%
         \advance\@itemnum by 1
         \expandafter \xdef\csname\@citelabel{##1}\endcsname{\number\@itemnum}%
         \ifx\biblabelprecontents\@undefined
            \let\biblabelprecontents = \hss
         \fi
         \ifx\biblabelpostcontents\@undefined
            \let\biblabelpostcontents = \relax
         \fi
         \@finishbibitem{##1}%
      }%
      \def\@finishbibitem##1{%
         \biblabelprint{\csname\@citelabel{##1}\endcsname}%
         \@writeaux{\string\@citedef{##1}{\csname\@citelabel{##1}\endcsname}}%
         \ignorespaces
      }%
      %
      %
      \let\em = \bblem
      \let\newblock = \bblnewblock
      \let\sc = \bblsc
      \frenchspacing
      \clubpenalty = 4000 \widowpenalty = 4000
      \tolerance = 10000 \hfuzz = .5pt
      \everypar = {\hangindent = \biblabelwidth
                      \advance\hangindent by \biblabelextraspace}%
      \bblrm
      \parskip = 1.5ex plus .5ex minus .5ex
      \biblabelextraspace = .5em
      \bblhook
      \input \bblfilebasename.bbl
   \endgroup
}%
%
%
\@innernewdimen\biblabelwidth
\@innernewdimen\biblabelextraspace
%
%
%
\def\biblabelprint#1{%
   \noindent
   \hbox to \biblabelwidth{%
      \biblabelprecontents
      \biblabelcontents{#1}%
      \biblabelpostcontents
   }%
   \kern\biblabelextraspace
}%
%
%
%
\def\biblabelcontents#1{{\bblrm [#1]}}%
%
%
\def\bblrm{\rm}%
%
%
\def\bblem{\it}%
%
%
\def\bblsc{\ifx\@scfont\@undefined
              \font\@scfont = cmcsc10
           \fi
           \@scfont
}%
%
%
\def\bblnewblock{\hskip .11em plus .33em minus .07em }%
%
%
\let\bblhook = \empty
%
%
%
\def\printcitestart{[}
\def\printcitefinish{]}
\def\printbetweencitations{, }
\def\printcitenote#1{, #1}
%
%
%
\let\citation = \@gobble
%
%
%
\@innernewcount\@numparams
%
%
\def\newcommand#1{%
   \def\@commandname{#1}%
   \@getoptionalarg\@continuenewcommand
}%
%
%
\def\@continuenewcommand{%
   \@numparams = \ifx\@optionalarg\empty 0\else\@optionalarg \fi \relax
   \@newcommand
}%
%
%
\def\@newcommand#1{%
   \def\@startdef{\expandafter\edef\@commandname}%
   \ifnum\@numparams=0
      \let\@paramdef = \empty
   \else
      \ifnum\@numparams>9
         \errmessage{\the\@numparams\space is too many parameters}%
      \else
         \ifnum\@numparams<0
            \errmessage{\the\@numparams\space is too few parameters}%
         \else
            \edef\@paramdef{%
               \ifcase\@numparams
                  \empty  No arguments.
               \or ####1%
               \or ####1####2%
               \or ####1####2####3%
               \or ####1####2####3####4%
               \or ####1####2####3####4####5%
               \or ####1####2####3####4####5####6%
               \or ####1####2####3####4####5####6####7%
               \or ####1####2####3####4####5####6####7####8%
               \or ####1####2####3####4####5####6####7####8####9%
               \fi
            }%
         \fi
      \fi
   \fi
   \expandafter\@startdef\@paramdef{#1}%
}%
%
%
%
%
\def\@readauxfile{%
   \if@auxfiledone \else 
      \global\@auxfiledonetrue
      \@testfileexistence{aux}%
      \if@fileexists
         \begingroup
            \endlinechar = -1
            \catcode`@ = 11
            \input \jobname.aux
         \endgroup
      \else
         \message{\@undefinedmessage}%
         \global\@citewarningfalse
      \fi
      \immediate\openout\@auxfile = \jobname.aux
   \fi
}%
%
%
\newif\if@auxfiledone
\ifx\noauxfile\@undefined \else \@auxfiledonetrue\fi
%
%
%
%
\@innernewwrite\@auxfile
\def\@writeaux#1{\ifx\noauxfile\@undefined \write\@auxfile{#1}\fi}%
%
%
%
\ifx\@undefinedmessage\@undefined
   \def\@undefinedmessage{No .aux file; I won't give you warnings about
                          undefined citations.}%
\fi
%
%
\@innernewif\if@citewarning
\ifx\noauxfile\@undefined \@citewarningtrue\fi
%
%
%
\catcode`@ = \@oldatcatcode

\def\pfeilso{\leavevmode
            \vrule width 1pt height9pt depth 0pt\relax
           \vrule width 1pt height8.7pt depth 0pt\relax
           \vrule width 1pt height8.3pt depth 0pt\relax
           \vrule width 1pt height8.0pt depth 0pt\relax
           \vrule width 1pt height7.7pt depth 0pt\relax
            \vrule width 1pt height7.3pt depth 0pt\relax
            \vrule width 1pt height7.0pt depth 0pt\relax
            \vrule width 1pt height6.7pt depth 0pt\relax
            \vrule width 1pt height6.3pt depth 0pt\relax
            \vrule width 1pt height6.0pt depth 0pt\relax
            \vrule width 1pt height5.7pt depth 0pt\relax
            \vrule width 1pt height5.3pt depth 0pt\relax
            \vrule width 1pt height5.0pt depth 0pt\relax
            \vrule width 1pt height4.7pt depth 0pt\relax
            \vrule width 1pt height4.3pt depth 0pt\relax
            \vrule width 1pt height4.0pt depth 0pt\relax
            \vrule width 1pt height3.7pt depth 0pt\relax
            \vrule width 1pt height3.3pt depth 0pt\relax
            \vrule width 1pt height3.0pt depth 0pt\relax
            \vrule width 1pt height2.7pt depth 0pt\relax
            \vrule width 1pt height2.3pt depth 0pt\relax
            \vrule width 1pt height2.0pt depth 0pt\relax
            \vrule width 1pt height1.7pt depth 0pt\relax
            \vrule width 1pt height1.3pt depth 0pt\relax
            \vrule width 1pt height1.0pt depth 0pt\relax
            \vrule width 1pt height0.7pt depth 0pt\relax
            \vrule width 1pt height0.3pt depth 0pt\relax}

\def\pfeilsw{ \leavevmode 
            \vrule width 1pt height0.3pt depth 0pt\relax
            \vrule width 1pt height0.7pt depth 0pt\relax
            \vrule width 1pt height1.0pt depth 0pt\relax
            \vrule width 1pt height1.3pt depth 0pt\relax
            \vrule width 1pt height1.7pt depth 0pt\relax
            \vrule width 1pt height2.0pt depth 0pt\relax
            \vrule width 1pt height2.3pt depth 0pt\relax
            \vrule width 1pt height2.7pt depth 0pt\relax
            \vrule width 1pt height3.0pt depth 0pt\relax
            \vrule width 1pt height3.3pt depth 0pt\relax
            \vrule width 1pt height3.7pt depth 0pt\relax
            \vrule width 1pt height4.0pt depth 0pt\relax
            \vrule width 1pt height4.3pt depth 0pt\relax
            \vrule width 1pt height4.7pt depth 0pt\relax
            \vrule width 1pt height5.0pt depth 0pt\relax
            \vrule width 1pt height5.3pt depth 0pt\relax
            \vrule width 1pt height5.7pt depth 0pt\relax
            \vrule width 1pt height6.0pt depth 0pt\relax
            \vrule width 1pt height6.3pt depth 0pt\relax
            \vrule width 1pt height6.7pt depth 0pt\relax
            \vrule width 1pt height7.0pt depth 0pt\relax
            \vrule width 1pt height7.3pt depth 0pt\relax
            \vrule width 1pt height7.7pt depth 0pt\relax
            \vrule width 1pt height8.0pt depth 0pt\relax
            \vrule width 1pt height8.3pt depth 0pt\relax
            \vrule width 1pt height8.7pt depth 0pt\relax
            \vrule width 1pt height9pt depth 0pt\relax
      }


\def\widestnumber#1#2{}

\def\citewarning#1{\ifx\shlhetal\relax 
    \else
    \par{#1}\par
    \fi
}

\def\rm{\fam0 \tenrm}

\def\fakesubhead#1\endsubhead{\bigskip\noindent{\bf#1}\par}



%
%
%

%

\font\textrsfs=rsfs10
\font\scriptrsfs=rsfs7
\font\scriptscriptrsfs=rsfs5

\newfam\rsfsfam
\textfont\rsfsfam=\textrsfs
\scriptfont\rsfsfam=\scriptrsfs
\scriptscriptfont\rsfsfam=\scriptscriptrsfs

\edef\oldcatcodeofat{\the\catcode`\@}
\catcode`\@11

\def\Cal@@#1{\noaccents@ \fam \rsfsfam #1}

\catcode`\@\oldcatcodeofat


\expandafter\ifx \csname margininit\endcsname \relax\else\margininit\fi

\long\def\red#1\endred{}
\long\def\green#1\endgreen{}
\long\def\blue#1\endblue{}

\def\endred{ \unmatched endred! }
\def\endgreen{ \unmatched endgreen! }
\def\endblue{ \unmatched endblue! }

\ifx\latexcolors\undefinedcs\def\latexcolors{}\fi

\def\emptycs{}
\def\evaluatelatexcolors{%
        \ifx\latexcolors\emptycs\else
        \expandafter\xxevaluate\latexcolors\xxfertig\evaluatelatexcolors\fi}
\def\xxevaluate#1,#2\xxfertig{\setupthiscolor{#1}%
        \def\latexcolors{#2}}

\font\smallfont=cmsl7
\def\rutgerscolor{\ifmmode\else\endgraf\fi\smallfont
\advance\leftskip0.5cm\relax}
\def\setupthiscolor#1{\edef\tmptmpcs{\noexpand\bgroup\noexpand\rutgerscolor
\noexpand\def\noexpand\currentcolor{#1}%
\noexpand}%
\expandafter\let\csname#1\endcsname\tmptmpcs
\def\tmptmpcs{\checkColorUnmatched{#1}\popthecolor}
\expandafter\let\csname end#1\endcsname\tmptmpcs}

\def\checkColorUnmatched#1{\def\expectcolor{#1}%
    \ifx\expectcolor\currentcolor   
    \else \edef\failhere{\noexpand\tryingToClose '\currentcolor' with end\expectcolor}\failhere\fi}

\def\currentcolor{???}

\def\popthecolor{\ifmmode\else\endgraf\fi\egroup}

\expandafter\def\csname#1\endcsname{}

\evaluatelatexcolors

 \let\outerhead\head
 \def\head{\innerhead}
 \let\innerhead\outerhead

 \let\outersubhead\subhead
 \def\subhead{\innersubhead}
 \let\innersubhead\outersubhead

 \let\outersubsubhead\subsubhead
 \def\subsubhead{\innersubsubhead}
 \let\innersubsubhead\outersubsubhead

 \def\proclaim{\innerproclaim}
 \let\innerproclaim\outerproclaim

 %
 %
 %
 %

\def\demo#1{\medskip\noindent{\it #1.\/}}
\def\enddemo{\smallskip}

\def\remark#1{\medskip\noindent{\it #1.\/}}
\def\endremark{\smallskip}

\pageheight{8.5truein}
\topmatter
\title{Spectra of monadic second order sentences} \endtitle
\author {Saharon Shelah \thanks {\null\newline
This research was supported by the Israel Science
Foundation. Publication 817 \null\newline
I would like to thank 
Alice Leonhardt for the beautiful typing. \null\newline
} \endthanks} \endauthor 

\affil{The Hebrew University of Jerusalem \\
Einstein Institute of Mathematics\\
Edmond J. Safra Campus, Givat Ram
Jerusalem 91904, Israel
 \medskip
Department of Mathematics \\
Hill Center-Busch Campus \\
Rutgers, The State University of New Jersey \\
110 Frelinghuysen Road \\
Piscataway, NJ 08854-8019  USA} \endaffil

\abstract  For a monadic sentence $\psi$ in the finite vocabulary
we show that the spectra, the set of cardinalities of
models of $\psi$ is almost periodic under reasonable conditions.  The
first is that every model is so called ``weakly $k$-decomposable".
The second is that we restrict ourselves to a nice class of models
constructed by some recursion.
\endabstract
\endtopmatter
\document

\newpage

\head {Introduction} \endhead  \resetall 
 \spuriousreset
\bigskip

This continues \cite{GuSh:536} and was announced there.
For a monadic second order sentence $\psi$ in the language
with one unary functions and unary predicates,
the spectra of the sentence
(i.e., the set $\{||M||: M $ a finite model of
$\psi \}$ is (see [GuSh: 536])
 periodic, but this fail badly
when we allow, e.g. two unary functions.
In the second section we characterize the family of finite 
structures which really behave like the unary function case, i.e., the
proof works.
\nl
In section one we assume that a monadic second order sentence
satisfies: every model is not indecomposable, 
i.e., has a non trivial  decoposition in a weak sense (see Definition
\scite{1.1}).
 We conclude that the specra is not arbitrary, 
 mainly - there are no big gaps in it
(from some point on). This is of
course considerably weaker  conclusion
than what we know for the languages with
only a unary function (under a much weaker assumption) or in \S2.
\nl
This work was done when Gurevich was writing \cite{GuSh:536}, but he
at first did not include an announcement in the version he circulated
insisting that I rewrote it to his satisfaction.
Meanwhile Fischer and Makowsky \cite{FiMw03} started
\cite{FiMw03} to work from the earlier version of \cite{GuSh:536}
continuing it in a different direction, using counting monadic logic
and dealing with tree and clique width of graphs (and of models).
It seems that Definition \scite{2.2} maybe a variant of ``clique width of
models"; see on this \cite{FiMw03}.
\nl
Clearly we can in \S1 use operations like \scite{2.2} instead of $M_1
\cup M_2$.
\bn
Note that in the definition of weakly $k$-decomposable we do not
require that the ``component" $M_1,M_2$ belongs to ${\frak K}$.  As it
was indirectly asked and to clarify Definition \scite{1.2}, we add:
\nl
\margintag{0.1}\ub{\stag{0.1} Example}:  The class ${\frak K}$ of finite incidence (or
edge) graphs is not weakly $k$-decomposable for every $k$.
\sn
Why?  Let $m$ be such that $m \ge k + \binom k2,m \ge 2k+2$ 
and it suffices to 
prove that for every
$n> m$, the statement $\circledast$ of Definition \scite{1.1}(1)
fail.  Let $G$ be the complete graph with set of nodes
$\{b_1,\dotsc,b_n\}$, so the incidence graph $G'$ has set of nodes $A =
\{b_1,\dotsc,b_n\} \cup \{c_{i,j}:1 \le i < j \le n\}$ and set of edges
$\{\{b_i,c_{i,j}\}:1 \le i < j \le n\} \cup \{\{b_j,c_{i,j}\}:1 \le
i < j \le n\}$.  So toward contradiction assume 
$A = A_1 \cup A_2,|A_1 \cap A_2| \le k$ and
$|A_1| \ge m,|A_2| \ge m$ with no edge (of $G'$) between $A_1 \backslash A_2$
and $A_2 \backslash A_1$.  Let $u = \{i:b_i \in A_1 \cap A_2$ or for
some $j,c_{i,j} \in A_1 \cap A_2$ or $c_{j,i} \in A_1 \cap A_2\}$.  
So $|u| \le 2k$.  Let $c_{\{i,j\}} = c_{\{j,i\}} = c_{i,j}$ when $1
\le i < j \le n$.
\bn
\ub{Case 1}:  For some $i_1,i_2$ we have $b_{i_1} \in A_1 \backslash
A_2,b_{i_2} \in A_2 \backslash A_1$.  

So for every $i \in \{1,\dotsc,n\} \backslash \{i_1,i_2\}$ for 
some $\ell \in \{1,2\},b_i \in A_\ell$, so
$\{b_i,c_{\{i,i_{3-\ell}\}}\}$ is an edge of $G'$ hence either $b_i \in
A_1 \cap A_2$ or $c_{i,i_{3-\ell}} \in A_1 \cap A_2$ hence (in both
possibilities) $i \in u$ so $n-2 \le |u|$ but $|u| \le 2k$ 
and $2k+2 \le m < n$, contradiction.
\bn
\ub{Case 2}:  Not Case 1.

So for some $\ell \in \{1,2\}$ we have $\{b_1,\dotsc,b_n\} \subseteq
A_\ell$,  hence $A_{3 - \ell} \backslash A_\ell \subseteq
\{c_{i,j}:b_i,b_j \in A_1 \cap A_2\}$, \wilog \, $\ell =1$ so 
$|A_2| \le |A_2 \cap A_1| +
|A_2 \backslash A_1| \le k + \binom k2 < m$, contradiction.
\hfill$\square_{\scite{0.1}}$\margincite{0.1}
\bigskip

\demo{\stag{0.2} Notation}  1) Let $n,m,\ell,k,i,j,d$ be natural
numbers.
\nl
2) $\tau$ is a vocabulary (i.e., set of predicates, individual
constants and function symbols, the last are not used here).
\nl
We thank Mor Doron for various helpful comments.
\enddemo
\newpage

\head {\S1 Weakly decomposable} \endhead  \resetall \sectno=1
 \spuriousreset
\bigskip

We can deal just with graphs just as this is traditional.  The
restriction to relational vocabulary is for simplifying our statements. 
\bigskip

\demo{\stag{1.0} Context}  1) Let $\tau$ be a finite relational
vocabulary, i.e., a finite set of predicates.
\nl
2) Let ${\frak K}^*_\tau$ be the class of $\tau$-models and recall
$\|M\|$ is the number of elements of $M \in {\frak K}^*_\tau,R^M$ is the
interpretation of the predicate $R \in \tau$. \nl
3) Let ${\frak K}$ denote a family of $\tau$-models 
closed under isomorphisms.
\enddemo
\bigskip

\definition{\stag{1.1} Definition}   1) We say that ${\frak K}$ is weakly
$k$-decomposable \ub{if}: for every $m$ there is $n$ such that
\mr
\item "{$\circledast_{k,m,n}$}"  if $M \in {\frak K},\|M\| \ge n$
\ub{then} we can find submodels $M_1,M_2$ (for graphs-induced 
subgraphs $G_1,G_2$) such that
{\roster
\itemitem{ $(a)$ }  $M_1 \cup M_2 = M$ i.e., $a \in M \Leftrightarrow
a \in M_1 \vee a \in M_2$ and $R^M = R^{M_1} \cup R^{M_2}$ for 
any $R \in \tau$ (for graphs: $G,G_1,G_2$ let $G = G_1 \cup G_2$ mean
that the set of nodes is the union
of the set of nodes of $G_1$ and of $G_2$, and the set of edges of $G$
is the union of the set of edges of $G_1$ and of $G_2$)
\sn
\itemitem{ $(b)$ }  $|M_1 \cap M_2| \le k$
\sn
\itemitem{ $(c)$ }  $\|M_\ell\| \ge m$ for $\ell = 1,2$.
\endroster}
\ermn
2) For a monadic second order sentence $\psi$ (in a vocabulary $\tau$)
we say that $\psi$ is weakly $k$-decomposable \ub{if} ${\frak K}^\tau_\psi$
is (see part (3)). \nl
3) For a vocabulary $\tau$ (as in \scite{1.0}) and sentence $\psi$ (in
this vocabulary) let ${\frak K}^\tau_\psi = \{M:M$ is a finite
$\tau$-model such that $M \models \psi\}$.  We may suppress $\tau$,
when clear from the context.
\enddefinition
\bigskip

\proclaim{\stag{1.2} Claim}   Assume
\mr
\item "{$(*)^{k^*}_\psi$}"  $\psi$ a monadic second order sentence, 
in the vocabulary $\tau$ such that ${\frak K} = {\frak K}^\tau_\psi$
is weakly $k^*$-decomposable
\ermn
\ub{then} ${\text{\rm Sp\/}}(\psi) = \{\|M\|:M \in {\frak K}\}$ 
satisfies for some $n^*$, that
\mr
\item "{$\circledast$}"  if $n_1 < n_2$ are successive members of
{\rm Sp}$(\psi)$ and $n^* < n_1$ then $n_2 < 2n_1$.
\endroster
\endproclaim
\bigskip

\demo{Proof}  Let $\psi$ have quantifier depth $\le d^*$. \nl
Let $m^*_1 > k^*$ be large enough such that
\mr
\item "{$\boxdot_1$}"  if $M_1 \in {\frak K}^*_\tau$ (yes ${\frak
K}^*_\tau$ and not ${\frak K}$), $\|M_1\| > k^*$ and
$a_1,\dotsc,a_k \in M_1$ and $k \le k^*$ \ub{then} 
there are $M_2 \in {\frak K}^*_\tau$ and $b_1,\dotsc,b_k \in M_2$ such that
(see \cite{GuSh:536} or \scite{2.6} below)
\nl
$\text{Th}^{d^*}(M_1,a_1,\dotsc,a_k) = 
\text{ Th}^{d^*}(M_2,b_1,\dotsc,b_k)$ \nl
and $k^* < \|M_2\| < m^*_1$.
\ermn
Let $m^*_2$ be such that the statement 
$\circledast_{k^*,m^*_1,m^*_2}$ from Definition \scite{1.1} holds (for
${\frak K}$).

Now assume that $n_1 < n_2$ are successive members of Sp$(\psi)$ and
$n_2 > m^*_2$.  Hence there is $M \in {\frak K}$ with exactly 
$n_2$ members.
So applying $\circledast_{k^*,m^*_1,m^*_2}$ to $M$ 
we can find $M_1,M_2$ as in Definition \scite{1.1} and let
$\{a_1,\dotsc,a_k\}$ list $M_1 \cap M_2$; so $k \le k^*$ and 
$\|M_1\|,\|M_2\| \ge m^*_1$.  Without loss of generality $\|M_1\| \le
\|M_2\|$, still $\|M_1\| \ge m^*_1$. \nl
By the choice of $m^*_1$ there is $(M'_1,b_1,\dotsc,b_k)$ such that
$k^* < \|M'_1\| < m^*_1$ and

$$
\text{Th}^{d^*}(M'_1,b_1,\dotsc,b_k) = 
\text{ Th}^{d^*}(M_1,a_1,\dotsc,a_2).
$$
\mn
Without loss of generality $\ell \in \{1,\dotsc,k\} \Rightarrow 
b_\ell = a_\ell$ and no member (for graphs - node) of 
$M'_1$ belongs to $M_2 \backslash \{a_1,\dotsc,a_k\}$.  Let
$M' = M'_1 \underset{\{a_1,\dots,a_k\}} {}\to + M_2$ be defined naturally
(set of elements of $M'  =$ union of set of elements of 
$M'_1$ and set of elements of $M_2$ and $R^M = R^{M'_1} 
\cup R^{M_2}$ for $R \in \tau$).

By the addition theorem (for local monadic theories, see
\scite{2.7}(c)) we have 
$M' \models \psi$, i.e., $M' \in {\frak K}$ and
\mr
\item "{$\boxdot$}"  $\frac 12 \|M\| \le \|M_2\| 
< \|M'\| = \|M'_1\| + \|M_2\| - k < m^*_1 + \|M_2\| - k \le
\|M_1\| + \|M_2\| - k = \|M\|$.
\ermn
That is $n_2/2 < \|M'\| < n_2$ but $M' \in {\frak K}$ so $\|M'\| \in
\text{ Sp}({\frak K})$ so there is $n' \in \text{ Sp}({\frak K})$ such
that $n_2/2 < n' < n_2$ so we are done.  \hfill$\square_{\scite{1.2}}$\margincite{1.2}
\enddemo
\bigskip

\demo{\stag{1.3} Conclusion}  If $\varphi$ is a second order monadic
sentence and $(*)^{k^*}_\varphi$ of \scite{1.2} holds and $\alpha$ is a
real $>0$ \ub{then} for some $n^*$ we have

$$
\boxtimes = \boxtimes_{\varphi,\alpha,n^*} \qquad \quad
n^* < n \in \text{ Sp}(\varphi) = (\exists m \in (\text{Sp}(\varphi))[n < m 
< (1 + \alpha)n].
$$
\enddemo
\bigskip

\demo{Proof}  Let $\Xi$ be the family of positive reals $\alpha$ such that
\mr
\item "{$\circledast_1$}"  for every monadic second order sentence
$\psi$ (for any vocabulary $\tau$ as in \scite{1.0}) such that
$(*)^{k^*}_\psi$ holds, the conclusion
$\boxtimes_{\varphi,\alpha,n^*}$ 
of \scite{1.3} holds for some $n^*$ (no harm in varying $k^*$, too).
\ermn
Note that allowing individual constants in $\tau$ is O.K. (either
allow them or code them by unary predicates); for a
vocabulary $\tau$ let $\tau^{+ k}$ be $\tau \cup
\{P_\ell:\ell=1,\dotsc,n\}$, where the $P_\ell$ are distinct unary
predicates not from $\tau$.
\nl
Clearly $0 < \beta < \alpha \and \beta \in \Xi \Rightarrow \alpha \in
\Xi$.  By Claim \scite{1.2} we have $1 \in \Xi$.

We shall now prove that
\mr
\item "{$\circledast_2$}"  if $\alpha \in \Xi \Rightarrow \alpha/2 \in
\Xi$. 
\ermn
This clearly suffices.  So let $\alpha,\tau,\psi,{\frak K} = 
{\frak K}^\tau_\psi$ be given.  Let $d$ be above the quantifier 
depth of $\psi$.  For $k \le k^*$ let

$$
\align
{\frak K}'_k = {\frak K}'_{\psi,k} = 
\{(M',P_1,\dotsc,P_k):&\text{for some } M \in {\frak K}
\text{ and } M_1,M_2 \text{ as in \scite{1.1}} \\
  &\text{in particular } |M_1 \cap M_2| \le k^* \text{ we have } M' = M_1 \\
  &\text{and } \{a_1,\dotsc,a_k\} \text{ lists } M_1 \cap M_2
\text{ and } P_\ell = \{a_\ell\}\}.
\endalign
$$
\mn
This is a class of $\tau^{+k}$ models.  Let
$\{\text{Th}^d(M',P_1,\dotsc,P_k):(M',P_1,\dotsc,P_k) \in {\frak
K}'_k\}$ be listed as $\bold t^k_1,\dotsc,\bold t^k_m$ and for $\ell \in
\{1,\dotsc,m\}$ let ${\frak K}'_{k,\ell} = {\frak K}'_{\psi,k,\ell} =
 \{(M',P_1,\dotsc,P_k) \in
{\frak K}'_k:\text{Th}^d(M',P_1,\dotsc,P_k) = \bold t^k_\ell\}$.
So $N \in {\frak K}'_{k,\ell} \Rightarrow |P^N_1| = \ldots = |P^N_k| = 1$.
\nl
It is not hard to see
\mr
\item "{$\circledast_3$}"  for some monadic second order sentence
$\psi_{k,\ell}$ of quantifier depth $d$ in the vocabulary $\tau^{+k},
{\frak K}'_{k,\ell}$ is the class of models of $\psi_{k,\ell}$, i.e.,
is ${\frak K}^{\tau^{+k}}_{\psi_{k,\ell}}$ for every relevant pair $(k,\ell)$.
\nl
[Why?  By direct checking ${\frak K}'_{k,\ell}$ is a class of
$\tau^{+k}$-models and there is such monadic sentences by the
definition of Th$^d$.]
\sn
\item "{$\circledast_4$}"  ${\frak K}^{\tau^{+k}}_{\psi_{k,\ell}}$ is
weakly $k^*$-decomposable.
\nl
[Why?  Clearly $\psi_{k,\ell}$ is a monadic sentence in the vocabulary
$\tau^{+k}$.
\nl
By the choice of $\bold t^k_\ell$ there are $M,M_1,M_2,a_1,\dotsc,a_k$
as in the definition of ${\frak K}'_k$.  Now expand $M_1$ to a
$\tau^{+k}$-model $M^*_1$ by $P^{M^*_1}_i = \{a_i\}$ and so
$\bold t^k_\ell = \text{ Th}^d(M^*_1)$.  We have to prove ``${\frak
K}'_{k,\ell}$ weakly decomposable", i.e., Definition \scite{1.2}.  So
let a number $m$ be given.  Let $m' = m + \|M_2\|$ and let $n$ be as
guaranteed for $m'$ by $(*)^{k^*}_\psi$ for ${\frak K}^\tau_\psi$.  
We shall show that $n$ is as
required for ${\frak K}^{\tau^{+k}}_{\psi_{k,\ell}}$.  Let
$(M'_1,P'_1,\dotsc,P'_k) \in {\frak K}'_{k,\ell},\|M'_1\| \ge n$ so
\wilog \, (i.e. by renaming) $P'_i = \{a_i\}$ (for
$i = 1,\dotsc,k$) and $M'_1 \backslash \{a_1,\dotsc,a_k\} \cap M_2 =
\emptyset$.  We can define $N$ such that $N,M'_1,M_2,a_1,\dotsc,a_k$
are as in \scite{1.2} so by the addition theorem (see \S2)
Th$^d(N,P'_1,\dotsc,P'_k) = \text{ Th}^d(M,P'_1,\dotsc,P'_k)$ so $N \models
\psi$.  As $\|N\| \ge \|M'_1\| \ge n$ by the choice of $n$ we can find
find $N_1,N_2,k' \le k^*$ and $\bar b = \langle
b_\ell:\ell=1,\dotsc,k' \rangle$ such that the tuple $(N,N_1,N_2,\bar
b)$ is as in \scite{1.2}, i.e., $N = N_1
\underset{\{b_1,\dotsc,b_{k'}\}} {}\to + N_2$
and $\|N_1\|,\|N_2\| \ge m'$.  Let $N'_\ell = N_\ell \restriction
(|M'_1|)$ and $\bar c = \langle c_\ell:\ell \ne 1,\dotsc,k''\rangle$
enumerate $\{b_\ell:\ell=1,\dotsc,k'\} \cap M'_1$.

Clearly $\|N'_\ell\| \ge \|N_\ell\| - \|M_2\| \ge m' - \|M_2\| = m$ by
the choice of $m'$ above.  So $(M'_1,N'_1,N'_2,\bar c)$ is as required
in the conclusion of \scite{1.2} for $\psi_{k,\ell}$.]
\sn
\item "{$\circledast_5$}"  $(*)^{k^*}_{\psi_{k,\ell}}$ holds
\nl
[why?  By $\circledast_3 + \circledast_4$.]
\ermn
Hence by the induction hypothesis for some $m^*$
\mr
\item "{$\circledast_6$}"   the conclusion of
$\boxtimes_{\psi_{k,\ell},\alpha,m^*}$ of \scite{1.3} holds (for any
relevant $k,\ell$)
\ermn
Let $n^*$ be as in \scite{1.1}(1) for $k^*,m^*,{\frak K} = 
{\frak K}^\tau_\psi$.
\nl
So it is enough to prove that $\boxtimes_{\psi,\alpha,n^*}$ holds.
Now for any $M \in {\frak K}$ with $\ge n^*$ elements there are
$M_1,M_2,k,a_1,\dotsc,a_k$ as in Definition \scite{1.1}(1) such that
$m^* \le \|M_1\| \le \|M_2\|$.  So for some $\ell,(M_1,P_1,\dotsc,P_k)
\in {\frak K}'_{k,\ell},P_i = \{a_i\}$ for $i=1,\dotsc,k$, so 
as we are assuming $\circledast_6$ clearly
\mr
\item "{$\circledast_7$}"   we can choose $(M'_1,P'_1,\dotsc,P'_k) \in
{\frak K}'_{k,\ell}$ such that
\endroster
\medskip

$$
\frac{\|M_1\|}{1 + \alpha} < \|M'_1\| < \|M_1\|.
$$
\mn
Without loss of generality (by renaming)

$$
P'_i = P_i,M'_1 \cap M_2 = \{a_1,\dotsc,a_k\}.
$$
\mn
We can define $M' \in {\frak K}^\tau_\psi$ as in the proof of
\scite{1.2} so with universe $|M'_1| \cup |M_2|$.  \nl
Hence

$$
\align
\circledast_8 \qquad \|M\| > \|M'\| &= \|M_2\| + \|M'_1\| - k > \|M_2\| + 
\frac{\|M_1\|}{1 + \alpha} - k \\
  &= \frac{1}{(1 + \alpha)}(\|M_2\| + \alpha\|M_2\| + \|M_1\| -k
- \alpha k) \\
  &= \frac{1}{1 + \alpha}(\|M\| + \alpha\|M_2\| - \alpha k) \\
  &\ge \frac{1}{1 + \alpha}(\|M\| + \alpha\|M\|/2) - \frac{\alpha k}{1 +
\alpha} = \frac{1 + \alpha/2}{1 + \alpha} \|M\| - \frac{\alpha k}{1 + \alpha}.
\endalign
$$
\mn
Let

$$
\beta =: \frac{1 +\alpha}{1 + \alpha/2} - 1 =
\frac{\alpha/2}{1 + \alpha/2} < \alpha/2.
$$
\mn
So by $\circledast_8$ we have

$$
\align
\|M\| > \frac{1 + \alpha/2}{1+\alpha} \|M\| - \frac{\alpha k}{1
+\alpha} &= \frac{1}{(1 + \beta)} \|M\| - \frac{\alpha k}{1 + \alpha} =
\frac{1}{1+ \alpha/2} \|M\| \\
  &+ (\frac{1}{1 + \beta} - \frac{1}{1 + \alpha/2}) \|M\| - 
 \frac{\alpha k}{1+k} \\
  &=\frac{1}{1 + \alpha/2} \|M\| + (\frac{\alpha/2-\beta}{(1 +
 \beta)(1 + \alpha/2)} \|M\| - \frac{\alpha k}{1 + \alpha}).
\endalign
$$
\mn
So we conclude: if conclusion \scite{1.3} 
holds for $\alpha > 0$ it holds for $\alpha/2$ 
provided that $(\alpha/2-\beta)\|M\| > \frac{\alpha k}{1 + \alpha}(1
+ \beta)(1 + \alpha/2)$, of course which holds if $\|M\|$ is large enough.
So we can prove by induction on $i$ that $\boxtimes$ holds for $\alpha \ge
\frac{1}{2^i}$.  \hfill$\square_{\scite{1.3}}$\margincite{1.3}
\enddemo
\newpage

\head {\S2 what the method of \cite{GuSh:536} gives} \endhead  \resetall \sectno=2
 \spuriousreset
\bn
\margintag{2.1}\ub{\stag{2.1} Discussion}:  The result above is interesting but leave
us unsatisfied.  For trees we get essentially sharp results.  Here the
spectra is not characterized.  We know that it is quite restricted
but, e.g. is it almost periodic?

The problem is that we do not see here a parallel to the operations
generating the class.
\bn
We may consider such classes:
\definition{\stag{2.2} Definition}  Let $\tau$ and $k^*$ be fixed, $\tau$
a finitary vocabulary with predicates only (coding function and 
individual constants by them if necessary)
and let ${\frak K}_{k^*} = {\frak K}^{c \ell}_{\tau,k^*}$ be 
the minimal family of 
$(M,a_1,\dotsc,a_k),M$ a finite $\tau$-model, 
$k \le k^*,a_\ell \in M$ such that
\mr
\item "{$(a)$}"  ${\frak K}_{k^*} = {\frak K}_{\tau,k^*}$ 
includes all the $(M,c_1,\dotsc,c_k), k \le k^*,c_\ell \in M$ with $M$ a
$\tau$-structure with $\le k^*$ elements
\sn
\item "{$(b)$}"  if $(M_\ell,a^\ell_1,\dotsc,a^\ell_{k_\ell}) 
\in {\frak K}_{k^*}$ for $\ell = 1,2$ and \nl
$x \in M_1 \wedge x \in M_2 \Rightarrow x \in
\{a^1_1,\dotsc,a^1_{k_1}\} \cap \{a^2_1,\dotsc,a^2_{k_2}\}$ \ub{then}
$(M,b_1,\dotsc,b_k) \in {\frak K}_{k^*}$ when:
{\roster
\itemitem{ $\circledast$ } $(i) \quad x$ an element of $M \Rightarrow
x$ an element (= node) of $M_1$ or of $M_2$
\sn
\itemitem{ ${{}}$ }  $(ii) \quad x$ an element of $M_\ell,x \notin
\{a^\ell_1,\dotsc,a^\ell_{k_\ell}\} \Rightarrow x$ an element of $M$
\sn
\itemitem{ ${{}}$ }  $(iii) \quad \{b_1,\dotsc,b_k\} \subseteq 
\{a^1_1,\dotsc,a^1_{k_1}\} \cup \{a^2_1 \ldots a^2_{k_2}\}$
\sn
\itemitem{ ${{}}$ }  $(iv) \quad$ if $R \in \tau$ is $m$-place
predicate, and $y_1,\dotsc,y_m \in M,z_1,\dotsc,z_m$ \nl

$\qquad \qquad \in M$  \ub{then}
$\langle y_1,\dotsc,y_m \rangle \in R^M \equiv 
\langle z_1,\dotsc,z_m \rangle \in R^M$ when:
\sn
\itemitem{ ${{}}$ }  $\quad \quad (\boxdot) \quad
(z_i=z_j) \equiv (y_i=y_j),(z_i=a^1_\ell) \equiv (y_i = a^1_\ell)$ \nl

\hskip60pt $(z_i=a^2_\ell) 
\equiv (y_i = a^2_\ell),(z_i \in M_\ell) \equiv (y_i \in M_\ell)$ 
and letting \nl

\hskip60pt $w_\ell = \{i:y_i \in M_\ell\}$ the quantifier
free type of $\langle y_i:i \in w_\ell \rangle$  \nl

\hskip60pt in $M_\ell$ is equal to the quantifier free type \nl

\hskip60pt of $\langle z_i:i \in w_\ell \rangle$ in $M_\ell$ for $\ell=1,2$.
\endroster}
\endroster
\enddefinition
\bigskip

\proclaim{\stag{2.3} Claim}  We can prove for ${\frak K}_{\tau,k^*}$
what we have proved for trees in \cite{GuSh:536}; 
including almost periodically of the spectrum for monadic sentences
(see \scite{2.7}(f)).
\endproclaim
\bigskip

\demo{Proof}  This is clause (f) of Claim \scite{2.7} proved below 
(as in \cite{GuSh:536}).
\enddemo
\bn
\margintag{2.4}\ub{\stag{2.4} Question}:  Is the class $K_{k^*}$ known?
Interesting? (see \S0)

Of course, e.g., the result on the spectrum is inherited by reducts.
After second thoughts I decide to add details (but naturally in the
style of \cite{Sh:42} rather than \cite{GuSh:536} as far as there is a
difference).
\bn
\margintag{2.5}\ub{\stag{2.5} Notation}:  $\tau$ is a vocabulary which has only
predicates and possibly individual constants.  If $R \in \tau$ is an
$m$-place predicate we write arity$(R) = m$; let arity$(\tau) =
\max(\{1\} \cup \{\text{arity}(R):R \in \tau\})$.

Let $\tau_m = \tau + \{P_0,\dotsc,P_{m-1}\}$ be the 
vocabulary $\tau$ when we add the (new and pairwise distinct)
predicates $P_0,\dotsc,P_{m-1}$ which below will be unary, similarly $\tau +
\{c_0,\dotsc,c_{k-1}\}$ for $c_\ell$ individual constants and
$\tau_{m,k} = \tau_m + \{\tau_0,\dotsc,c_{k-1}\}$.
\bigskip

\definition{\stag{2.6} Definition}  1) Let $\tau$ be a finite
vocabulary $\tau$ consisting of predicates only; $P_0,P_1,\ldots$ be unary
predicates $\notin \tau$; (for notational simplicity).  For a
$\tau$-model $M$ and sequence $\bar{\Cal U}^m = \langle {\Cal
U}_\ell:\ell < m \rangle$ of subsets of $M$ we define
Th$^n(M,\bar{\Cal U}^m)$ by induction on $n$:
\mr
\item "{$(a)$}"  $n=0$ it is the set of sentences $\psi = (\exists
x_0,\dotsc,x_{r-1}) \wedge \Phi$ where $r \le$ arity$(\tau)+1,\Phi$ a 
set of basic formulas of $\tau + \{P_0,\dotsc,P_{m-1}\}$ such that 
$(M,\bar{\Cal U}^m) \models \psi$
\sn
\item "{$(b)$}"  Th$^{n+1}(M,\bar{\Cal U}^m) =
\{\text{Th}^n(M,\bar{\Cal U}^m \char 94 \langle {\Cal
U}_m\rangle):{\Cal U}_m \subseteq M\}$.
\ermn
2) TH$^{n,m}(\tau)$ is the set of formally possible Th$^n(M,\langle
{\Cal U}_\ell:\ell < m \rangle)$; see below; if $m=0$ we may omit $m$.
\nl
3) For a class ${\frak K}$ of models let
TH$^{n,m}({\frak K})$ be $\{\text{Th}^n(M,{\Cal U}_0,\dotsc,
{\Cal U}_{m-1}):M \in {\frak K},{\Cal U}_0,\dotsc,\subseteq M\}$; if
$m=0$ we may omit it.
\nl
4) Let ${\frak K}_{\tau,k^*,k}$ be the set of models
$(M,c_1,\dotsc,c_k) \in {\frak K}_{\tau,k^*},M$ a $\tau$-structure.
\enddefinition
\bigskip

\proclaim{\stag{2.7} Claim}  Let $\tau$ have predicates only, define
$\tau_k = \tau + \{c_0,\dotsc,c_{k-1}\}$ and let 
$k^* \ge { \text{\rm arity\/}}(\tau)$ and $n$ be given.

We can compute the following (from $\tau,k^*,n^*,m^*$)
\mr
\item "{$(a)$}"  for $n \le n^*,m \le m^* - n,k \le k^*$, {\rm
TH}$^n(\tau_{m,k})$ which, the set of formally possible 
{\rm Th}$^{n}(M,\langle {\Cal U}_0,\dotsc,{\Cal U}_{m-1}\rangle,
a_0,\dotsc,a_k),M$ a $\tau$-model 
\sn
\item "{$(b)$}"  the set ${\Cal S}_{\tau_m,k^*,k_1,k_2,k}$ of schemes
defined implicitly in Definition \scite{2.2}
\sn
\item "{$(c)$}"  we can compute the functions \nl
$F =: F^n_{\tau_m,k_1,k_2,k}:{\text{\rm TH\/}}^n(\tau_{m,k_1}) 
\times { \text{\rm TH\/}}^n(\tau_{m,k_2}) \times {\Cal S}_{\tau_m,
k^*,k_1,k_2,k} \rightarrow { \text{\rm TH\/}}^n(\tau_{m,k})$ such that
(where $M_1,M_2,M$  are $\tau_m$-models) if
$(M_\ell,a^\ell_0,\dotsc,a^\ell_{k_\ell-1})$ for $\ell=1,2$ and
$(M,a_0,\dotsc,a_{k-1})$ are as in Definition \scite{2.2} for the
scheme $\bold s$ and $t_\ell = { \text{\rm
Th\/}}^n(M_\ell,a^\ell_0,\dotsc,a^\ell_{k_\ell -1})$ for $\ell=1,2$
and $t = { \text{\rm Th\/}}^n(M,a_0,\dotsc,a_{k-1})$ \ub{then} $t =
F^n(t_1,t_2,\bold s)$, in particular the representative models does
not matter
\sn
\item "{$(d)$}"  we can compute $T^n_{\tau_m,k^*,k} = \{{\text{\rm
Th\/}}^n(M,c_0,\dotsc,c_{k-1}):M$ is a $\tau_m$-model with $\le k^*$
elements$\}$
\sn
\item "{$(e)$}"  the sequence $\langle{\text{\rm
TH\/}}^n({\frak K}_{\tau,k^*,k}):k \le k^* \rangle$ can be computed
(could use $\tau_m$)
\sn  
\item "{$(f)$}"  for each $t \in { \text{\rm TH\/}}^n({\frak
K}_{\tau,k^*,k})$, the set

$$
{\text{\rm Sp\/}}_t = \{\|M\|:(M,a_0,\dotsc,a_{k-1}) \in
{\frak K}_{\tau,k^*,k} \text{ and } {\text{\rm
Th\/}}^n(M,a_0,\dotsc,a_{k-1})=t\}
$$
is eventually periodic; i.e. for some $m_1,m_2 \in M$ we have: if
$\ell_1,\ell_2 > m_1$ are equal modulo $m_2$ \ub{then} $\ell_1 \in {
\text{\rm Sp\/}}_t \equiv \ell_2 \in { \text{\rm Sp\/}}_t$; 
\endroster
\endproclaim
\bigskip

\demo{Proof}  Straight by now, e.g.,
\sn
\ub{Clause (a)}:  For $n = 0$ trivial and for $n+1$ we can let

$$
{\text{\rm TH\/}}^{n+1}(\tau_{m,k}) = 
{\Cal P}({\text{\rm TH\/}}^n(\tau_{m+1,k})).
$$
\sn
\ub{Clause (c)}:  This is the addition theorem proved by induction on
$n$.
\sn
\ub{Clause (e)}: We start with the sequence $\langle T^n_{\tau,k^*,k}:k \le
k^* \rangle$ and close it under the operations
$F^n_{\tau,k^*,k_1,k_2}(-,-,\bold s),\bold s \in 
{\Cal S}_{\tau,k^*,k_1,k_2,k}$ .
That is we define $\langle T^i_k:k \le k^* \rangle$ by induction on $i
\le \dsize \sum_{k<k^*}|{\text{\rm TH\/}}^n(\tau_k)|$ as follows:
\mr
\item "{$(\alpha)$}"   $T^0_k = T^n_{\tau,k^*,k}$ from clause (d)
\sn  
\item "{$(\beta)$}"  $T^{i+1}_k = T^i_k \cup
\{F^n_{\tau,k^*,k_1,k_2}(t_1,t_2,\bold s),t_1 \in T^i_{k_1},t_2 \in
T^i_{k_2},\bold s \in {\Cal S}^n_{\tau,k^*,k_1,k_2,k}\}$. 
\ermn
By cardinality considerations we know that for the last $i$ 
we have ${\text{\rm TH\/}}^n(K_{\tau,k^*,k}) =
T^i_k$ (for $k \le k^*$) as
\mr
\item "{$(\gamma)$}"  for each $k,\emptyset \subseteq T^0_k \subseteq
T^1_k \subseteq \ldots \subseteq T^i_k \subseteq \ldots$ and the number
of $i$ for which $T^{k_i}_k \ne T^{i+1}_k$ is $\le
|\text{TH}^n(\tau_k)|$ and
\sn
\item "{$(\delta)$}"  if $T^i_k = T^{i+1}_k$ for every $k \le k^*$
then $j \ge i \wedge k \le k^* \Rightarrow T^j_k = T^i_k$.
\ermn
\ub{Clause (f)}:  By the proof above and observation \scite{2.8} below.
\hfill$\square_{\scite{2.7}}$\margincite{2.7}
\enddemo
\bigskip

\demo{\stag{2.8} Observation}  Assume that
\mr
\item "{$(a)$}"  $m \in \Bbb N,\bold W$ is a finite set of quadruples
$\{\langle \ell_1,\ell_2,\ell_3,j \rangle:\ell_1,\ell_1,
\ell_3 < m$ and $j \in \Bbb N\}$
\sn
\item "{$(b)$}"  $N^i_\ell \, (\ell < m,i \in \Bbb N)$ is a finite set
of natural numbers
\sn
\item "{$(c)$}"  $N^{i+1}_\ell = N^i_\ell \cup \{n_1 + n_2-j$: for
some $(\ell_1,\ell_2,\ell_3,j) \in \bold W$ we have $n_1 \in
N^i_{\ell_1},n_2 \in N^i_{\ell_2},\ell_3=\ell\}$.
\ermn
\ub{Then} each set $N_\ell =: \cup\{N^i_\ell:i$ a natural number$\}$ is 
almost periodic.  
\enddemo
\bigskip

\demo{Proof}  Let $i < m$, clearly $n \in N_i$ \ub{iff} we can
find a witness $({\Cal T},\bar \ell,\bar n,\bar w)$ which means
\mr
\item "{$\circledast$}"  $(a) \quad {\Cal T}$ is a finite set of
sequences of zeroes and ones closed under \nl

\hskip30pt initial segments
\sn
\item "{${{}}$}"  $(b) \quad \bar \ell = \langle \ell_\eta:\eta \in
{\Cal T} \rangle$ such that $\ell_{<>} = i$
\sn
\item "{${{}}$}"  $(c) \quad \bar n = \langle n_\eta:\eta \in {\Cal T}
\rangle$ such that $n_{<>} = n$
\sn
\item "{${{}}$}"  $(d) \quad \bar w = \langle w_\eta:\eta \in {\Cal T}$
not maximal$\rangle, w_\eta \in \bold W$
\sn
\item "{${{}}$}"  $(e) \quad$ if $\eta \in {\Cal T}$ is
$\triangleleft$-maximal then $n_\eta \in N^0_{\ell_\eta}$
\sn
\item "{${{}}$}"  $(f) \quad$ if $\eta \in {\Cal T}$ is not
$\triangleleft$-maximal then $\eta_0 = \eta \char 94 \langle 0 \rangle
\in {\Cal T}$ and $\eta_1 = \eta \char 94 \langle 1 \rangle \in {\Cal T}$ 
\nl

\hskip23pt and we have
$$
w_\eta = (\ell_{\eta \char 94 \langle 0 \rangle},
\ell_{\eta \char 94 \langle 1 \rangle},\ell_\eta,
n_{\eta \char 94 \langle 0 \rangle} + n_{\eta \char 94 \langle 1
\rangle} - n_\eta)
$$
\sn
\item "{${{}}$}"  $(g) \quad$ (follows): $n_\eta\in N_{\ell_\eta}$ for
$\eta \in {\Cal T}$.
\ermn
Now let $n^*_1 =: 2^m \times n^*_0$ where $n^*_0 =: 
\max(\dbcu_\ell N^0_\ell)$ and let $n^* =: n^*_1!$.
Assume
\mr
\item "{$(*)_0$}"  $n^* < n \in N_\ell$ and 
there is no $n'$ such that $n^* < n' \in N_\ell,n' < n$ 
and $n=n'$ mod $n^*$.
\ermn
Choose $({\Cal T},\bar \ell,\bar n,\bar w)$ as above such that 
$(\ell_{<>},n_{<>}) = (\ell,n)$ and $|{\Cal T}|$
minimal.  
Let

$$
{\Cal U} = \{\nu \in {\Cal T}:\text{if } \nu \triangleleft \rho 
\in {\Cal T} \text{ then } n_\nu < n_\rho\}
$$
\mn
and let

$$
t(\nu) = \text{ Max}\{|\{m:\nu \trianglelefteq \rho \restriction m
\text{ and } \rho \restriction m \in {\Cal U}\}|:\nu \triangleleft \rho
\in {\Cal T}\}.
$$
\mn
Now
\mr
\item "{$(*)_1$}"  $\nu_0 \trianglelefteq \nu_1 \in {\Cal T}
\Rightarrow t(\nu_0) \ge t(\nu_1)$ and if $\nu \in {\Cal U}$ is not
$\triangleleft$-maximal then $t(\nu) = \text{ Min}\{t(\nu {}^\frown \langle j
\rangle)+1:\nu {}^\frown \langle j \rangle \in {\Cal T}\}$.
\nl
[Why?  Look at the definitions.]
\sn
\item "{$(*)_2$}"  if $\nu \in {\Cal T}$ then $n_\nu \le 2^{t(\nu)}
\times n^*_0$
\nl
[Why?  We prove this by induction on $t(\nu)$ and then on $|\{\rho:\nu
\trianglelefteq \rho \in {\Cal T}\}|$.  If $\nu$ is
$\triangleleft$-maximal in ${\Cal T}$ necessarily $n_\nu \in
N^0_{\ell_\nu}$ hence $n_\nu \le n^*_0$ so the conclusion of $(*)_2$
in this case is trivial.  If $\nu \in {\Cal U},\nu$ is not
$\triangleleft$-maximal in ${\Cal T}$ then $n_\nu \le 2 
\times \text{ Max}\{n_{\nu {}^\frown <j>}:\nu {}^\frown <j> \in 
{\Cal T},j=0,1\}$ and clearly
$t(\nu) = \text{ Max}\{t(\nu {}^\frown <j>)+1:j=0,1\}$, so we can check easily.
Lastly, if $\nu \in {\Cal T} \backslash {\Cal U}$ we can find 
$\nu_1$ such that $n_\nu \le n_{\nu_1},\nu \triangleleft \nu_1 \in {\Cal T}$ so
$t(\nu) \ge t(\nu_1)$ by $(*)_1$ and so $n_\nu \le n_{\nu_1} \le
2^{t(\nu_1)} \times n^*_0 \le 2^{t(\nu)} \times n^*_0$ as required.]
\sn
\item "{$(*)_3$}"  there are $\nu_0 \triangleleft 
\nu_1$ from ${\Cal U}$ such that
$\ell_{\nu_0} = \ell_{\nu_1}$.
\nl
[Why?  As $n_{<>} = n > 2^m \times n^*_0$, by $(*)_2$ we know that
$t(<>) > m$, hence we can find $\nu_0 \triangleleft \nu_1
\triangleleft \ldots \triangleleft \nu_m$ all from ${\Cal U}$.  As
$\ell_{\nu_0},\dotsc,\ell_{\nu_m} < m$ clearly for some $j(0) <
j(1) \le m$ we have $\ell_{\nu_{j(0)}} = \ell_{\nu_{j(1)}}$.]
\sn
\item "{$(*)_4$}"  if $\nu_0 \in {\Cal U}$ is $\triangleleft$-maximal such that
for some $\nu_1,\nu_0 \triangleleft \nu_1 \in {\Cal U} \and 
\ell_{\nu_0} = \ell_{\nu_1}$ then $t(\nu_0) \le m$.
\nl
[Why?  Look at the definition of $t(\nu_1)$, i.e., repeat the proof of
$(*)_3$.]
\sn
\item "{$(*)_5$}"  for some $\nu_0 \triangleleft \nu_1$ from 
${\Cal U} \subseteq {\Cal T},\ell_{\nu_0} = 
\ell_{\nu_1}$ and $n_{\nu_0} < n_{\nu_1}$ and $n_{\nu_0} \le n^*_1$
\nl
[Why?  By $(*)_3$ we can find $\nu_0 \triangleleft \nu_1$ from ${\Cal
U}$ with $\ell_{\nu_0} = \ell_{\nu_1}$ so without loss of generality $\nu_1$ is
$\triangleleft$-maximal, hence by $(*)_4$ we know 
$t(\nu_0) \le m$, so by $(*)_2$,
$n_{\nu_0} \le n^*_1$.]
\ermn
Now we can take a copy of $A = \{\rho \in {\Cal T}:\nu_0 \trianglelefteq 
\rho$ but $\neg \nu_1 \triangleleft \rho\}$
and insert it just before $\nu_1$ any number of times, hence
$n+i(n_{\eta_0} - n_{\eta_1}) \in N_\ell$ for any $i$.  As $n_{\nu_0}
- n_{\nu_1} \le n_{\eta_1} < n^*_1$ we are done.  (We can compute a
bound when this starts.  That is omitting $A$ we get
$n-(n_{\eta_0}-n_{\eta_1}) \in N_\ell$ hence by the assumption on $n$
in $(*)_0,n =: n-(n_{\eta_0}-n_{\eta_1}) \le n^*$ so $n \le
n^*+(n_{\eta_0}-n_{\eta_1}) \le n^* + n^*_1$.)  \hfill$\square_{\scite{2.8}}$\margincite{2.8}
\enddemo
\bigskip

\remark{Remark}  Of course, also in \S1 we can use sums as in \scite{2.2}.
\endremark
\bigskip

\nocite{ignore-this-bibtex-warning} 
\newpage
    
REFERENCES.  
\bibliographystyle{lit-plain}
\bibliography{lista,listb,listx,listf,liste}

\enddocument